\documentclass{amsart}

\usepackage{pifont}
\usepackage{amsfonts}

\usepackage{amscd,amssymb,amsmath,graphicx,verbatim,mathrsfs}
\usepackage[TS1,OT1,T1]{fontenc}
\usepackage{extarrows}
\usepackage[all]{xy}
\usepackage{dcpic,mathtools,pictexwd}
\usepackage{rotating}
\theoremstyle{remark}
\usepackage{indentfirst}
\usepackage{amsthm} 
\usepackage{tikz}
\usetikzlibrary{angles,quotes}
\usetikzlibrary{calc}
\usepackage{amsmath}

\usepackage{color}
\usepackage[colorlinks=true, allcolors=blue]{hyperref}
\usepackage{cleveref}
\usepackage{soul}
\usepackage{enumitem}
\usepackage{tikz-cd}

\allowdisplaybreaks

\theoremstyle{plain}

\newtheorem{thm}{Theorem}[section]
\newtheorem{lem}[thm]{Lemma}

\newtheorem{prop}[thm]{Proposition}

\newtheorem{ex}[thm]{Example}

\theoremstyle{definition}

\newtheorem{defn}[thm]{Definition}

\newtheorem{rmk}[thm]{Remark}

\newtheorem{constr}[thm]{Construction}

\usepackage{setspace}

\usepackage{geometry}
\theoremstyle{plain}
\newcounter{thmintroctr}

\theoremstyle{definition}
\newcounter{goalintroctr}

\numberwithin{equation}{section}

\begin{document}

\linespread{1.2}

\title{Higher-order jet and its application to the Novikov conjecture}

\author{Jintao Deng}
\address{Department of mathematics, Center of Noncommutative Geometry, Shanghai University of Finance and Economics}
\email{jintaodeng@sufe.edu.cn}
\thanks{The first author is partially supported by Shanghai Pujiang Project, No. 25PJA033.}

\author{Geng Tian}
\address{Department of Mathematics, Liaoning University, Liaoning, China 110036}
\email{gengtian.ncg@gmail.com}

\thanks{The second author was supported in part by Grant for Excellent Young Scholars in Tianyuan Mathematics.}


\begin{abstract}
In this paper, we apply higher-order tangent bundles and higher-order derivative maps to study countable subgroups of the deffeomorphism group ${\rm Diff}(M)$ of a closed smooth manifold $M$. Using this framework, we introduce the concept of higher order $\mu$-discreteness for subgroup $\Gamma < {\rm Diff}(M)$.

As applications, we show that both the strong Novikov conjecture and the coarse Baum--Connes conjecture hold for countable groups with higher order $\mu$-discreteness.  
\end{abstract}

\maketitle
	
\tableofcontents

\section{Introduction}

The Novikov conjecture asserts that the higher signatures of closed
oriented smooth manifolds are invariant under orientation-preserving
homotopy equivalences. 
Noncommutative geometry provides an approach to the study of the
Novikov conjecture via higher index theory.
More precisely, for a countable discrete group $\Gamma$, the Novikov
conjecture for  closed oriented smooth manifolds  with fundamental
group $\Gamma$ follows from the rational strong Novikov conjecture
for $\Gamma$. The latter asserts that the rational
Kasparov--Baum--Connes assembly map
$$
K_*(B\Gamma)\otimes_{\mathbb Z}\mathbb Q
\longrightarrow
K_*\bigl(C_r^*(\Gamma)\bigr)\otimes_{\mathbb Z}\mathbb Q
$$
is injective, where $B\Gamma$ is the classifying space for $\Gamma$,
$K_*(B\Gamma)$ denotes its $K$-homology with compact supports,
and $C_r^*(\Gamma)$ is the reduced group $C^*$-algebra.
Although the Novikov conjecture has been established for many classes
of groups, it remains open in general for a natural and important class:
countable subgroups of diffeomorphism groups of smooth manifolds.

For a closed smooth manifold $M$, Gong, Wu and Yu
\cite{GWY} proved the rational strong Novikov conjecture
for countable subgroups $\Gamma<\operatorname{Diff}(M)$ that preserve
a smooth density $\mu$ 
(a special type of Borel measure on $M$)  and whose
natural isometric action on the $L^2$-direct integral
$$
L^2_{\textsl{g}_1}
\bigl(M,\mu,\operatorname{Riem}(M)\bigr)
$$
is metrically proper. 
Here  $\textsl{g}_1$  is a fixed 
Riemannian metric on $M$ and 
$\operatorname{Riem}(M)\to M$  is the 
Riemannian metric bundle of $M$. 
Their notion of \emph{geometric discreteness}
characterizes the metric properness of this action.

Subsequently, Gong, Wu, Xie and Yu \cite{GongWuXieYu} generalized the above result by removing the measure-preserving assumption. 
They proved the rational strong Novikov conjecture for countable subgroups $\Gamma<\operatorname{Diff}(M)$ for which the canonical orbit map
$$
\Gamma\longrightarrow
L^2_{\textsl{g}_1}
\bigl(M,\mu,\operatorname{Riem}(M)\bigr),
\qquad
\gamma\longmapsto\bigl[\gamma_* \textsl{g}_1\bigr]
\qquad \qquad (*)
$$
is a coarse embedding for some Borel probability measure $\mu$ on $M$,
where $\textsl{g}_1$ is a fixed continuous  Riemannian metric. The difficulty of understanding the large scale geometry of the group $\Gamma$ arises from the fact that the space 
$L^2_{\textsl{g}_1}
\bigl(M,\mu,\operatorname{Riem}(M)\bigr)$ is just a space of $L^2$-sections from the space $M$ to the fiber the Hilbert-Hadamard space they considered, not a Hilbert space. 

The notion of \emph{$\mu$-discreteness} introduced in \cite{GongWuXieYu} characterizes the coarse embedding property of this orbit map.
The orbit map $(\ast)$ derives from the natural isometric group action of  $\Gamma$ 
on the continuous field
$$
\left\{
L^2_{\textsl{g}_1}
\bigl(M,\nu,\operatorname{Riem}(M)\bigr)
\right\}_{\nu\in\operatorname{Prob}(M)},
$$
where $\operatorname{Prob}(M)$  is  the space of Borel probability measures on $M$, equipped with the weak$^*$-topology. Using the group action on this continuous field, they proved the
rational strong Novikov conjecture via a $KK$-theoretic argument involving the associated $C^*$-algebras. Their method is therefore limited to coarse embeddings of orbit maps.

Our approach in this paper overcomes this limitation and applies to general coarse
embeddings, without requiring them to be orbit maps.
Recall that  a map  $h:(X,d_X)\to(Y,d_Y)$ between metric spaces is a
\emph{coarse embedding}  \cite{Gromov1993}  
if there exist nondecreasing functions
$\rho_-,\rho_+:[0,\infty)\to[0,\infty)$, with
$\rho_-(t)\to\infty$ as $t\to\infty$, such that
$$
\rho_-\bigl(d_X(x,y)\bigr)
\leq d_Y\bigl(h(x),h(y)\bigr)
\leq \rho_+\bigl(d_X(x,y)\bigr),
\qquad  x,y\in X.
$$

In this paper, we study subgroups of diffeomorphism groups  using higher-order tangent
bundles and   higher-order derivatives of these diffeomorphisms.
In local coordinates, the ordinary derivative records only the linear term of the Taylor expansion of a diffeomorphism. However, the $r$-th order derivative encodes all partial derivatives of orders $1,\ldots,r$, thereby retaining the nonlinear terms through degree $r$. This additional information allows us to detect a strictly larger class of subgroups of $\operatorname{Diff}(M)$, compared with the $\mu$-discreteness of Gong--Wu--Xie--Yu \cite{GongWuXieYu}.

More precisely, for a compact smooth manifold $M$ equipped with a finite Borel measure
$\mu$, we introduce higher-order $\mu$-discreteness and higher-order orbit maps.
As in the first-order case, a countable subgroup $\Gamma<\operatorname{Diff}(M)$ is higher-order $\mu$-discrete if and only if the corresponding orbit map is a coarse embedding.
More generally, our approach extends beyond orbit maps:  the Novikov conjecture holds for every countable group $\Gamma$ that admits a coarse embedding into an $L^p$-direct integral
$$
L^p_{\textsl{g}_r}
\bigl(M,\mu,\operatorname{Riem}^{(r)}(M)\bigr)
$$
for some $p\in[1,\infty)$ and some integer $r\geq1$. Here, $\operatorname{Riem}^{(r)}(M)\to M$ denotes the $r$-th order Riemannian metric bundle of $M$. Restricting this result to the
higher-order orbit maps yields the Novikov conjecture for every
higher-order $\mu$-discrete subgroup of
$\operatorname{Diff}(M)$.
Our proofs adopt a direct approach based on coarse embeddings into
$L^{p'}$-spaces for suitable $p'\in[1,\infty)$, bypassing the
continuous field and $C^*$-algebraic constructions.
Hence the coarse Baum--Connes conjecture also holds for
every discrete metric space with bounded geometry  that admits a
coarse embedding into one of these $L^p$-direct integrals.
We emphasize  that this approach  applies to general $L^p$-direct integrals, which need
not be Hilbert--Hadamard spaces.

We emphasize that, even for $p=2$, the higher-order method detects
a strictly larger class of subgroups of diffeomorphism groups
than its first-order counterpart.
In the final section, we exhibit countable subgroups
$\Gamma<\operatorname{Diff}(M)$ for which the first-order orbit map
considered by Gong, Wu, Xie and Yu \cite{GongWuXieYu} fails to be
a coarse embedding for every Borel probability measure $\mu$ on $M$
and every continuous  Riemannian metric, 
whereas a
higher-order orbit map is a coarse embedding for a suitable choice
of measure. 
These examples include both amenable and non-amenable
groups and lie beyond the scope of their first-order criterion.

We  establish these results in the more general setting of
\emph{snowflake metric bundles}, whose typical fibers admit snowflake
embeddings into a separable Hilbert space. 
By the work of Lang, Schlichenmaier, Naor and Silberman \cite{LangSchlichenmaier,NaorSilberman2011}, we know that  besides the Hilbert space itself, the fibers of snowflake metric bundles can be homogeneous Hadamard manifolds, pinched negatively curved Hadamard manifolds, doubling spaces, separable metric trees, and separable Euclidean buildings, as well as finite products and finite unions of these spaces. 
In particular, the higher-order Riemannian metric bundles
$\operatorname{Riem}^{(r)}(M)\to M$ provide natural examples of
snowflake metric bundles.
For a general snowflake metric bundle, we prove that the Novikov
conjecture (respectively, the coarse Baum--Connes conjecture) holds
for every countable group (respectively, discrete metric space with
bounded geometry) that admits a coarse embedding into one of its
$L^p$-direct integrals.

The following is our main result.

\begin{thm}\label{main-theorem}
Let $M$ be a compact smooth manifold, and $\Gamma<\operatorname{Diff}(M)$ a higher-order $\mu$-discrete subgroup. Then the strong Novikov conjecture holds for $\Gamma$, i.e.
the Baum-Connes assembly map
$$
K_*(\underline{E}\Gamma)
\longrightarrow
K_*\bigl(C_r^*(\Gamma)\bigr)
$$
is injective.

In addition, if the group $\Gamma$ is equipped with proper left-invariant metric, then the coarse Baum--Connes conjecture holds for $|\Gamma|$.

\end{thm}

The paper is organized as follows. In Section~2, 
we introduce higher-order $\mu$-discreteness, and prove
the orbit-map characterization.
In Section~3,  we develop the theory
of snowflake metric bundles and prove Theorem~\ref{main-theorem}. 
Finally, in Section~4, we present
 examples for which the first-order criterion
fails but the higher-order criterion applies.

\emph{\bf AI disclosure.}  ChatGPT 6.0 Astra was used solely for language editing and grammatical
refinement during the preparation of this manuscript. All mathematical proofs and
core results were developed by the authors.

\section{Higher-order $\mu$-discreteness and higher-order orbit map}

In this section, we first introduce the concept of higher-order $\mu$-discreteness for countable subgroups of $\operatorname{Diff}(M)$, where
$M$ is a compact smooth manifold equipped with a finite
Borel measure $\mu$. This notion is defined in terms of the distortion of fiber inner products on higher-order tangent bundles under the
higher-order derivative maps induced by diffeomorphisms.
In the case when $r=1$ and $p=2$, our definition recovers the notion of
$\mu$-discreteness introduced by Gong--Wu--Xie--Yu
\cite{GongWuXieYu}.

We then construct the associated higher-order orbit map, and
show that its restriction to a countable subgroup
$\Gamma<\operatorname{Diff}(M)$, equipped with a proper
left-invariant metric, is a coarse embedding if and only if
$\Gamma$ is higher-order $\mu$-discrete.

\subsection{Higher-order derivative maps}

The ordinary derivative map between tangent bundles captures only the first-order linearization of a smooth map: in local coordinates, it retains the linear term of the Taylor expansion and discards all quadratic and higher-order terms. In contrast, the $r$-th order derivative map retains the full Taylor data through order $r$. More precisely, its value at a point is determined by all partial derivatives of orders $1,\ldots,r$, and therefore encodes the entire Taylor expansion of a smooth map truncated at degree $r$. Thus, at the first-order level one sees only infinitesimal linear behavior, whereas at higher-order levels one also sees the successive nonlinear terms of the Taylor expansion. 
This additional information allows us to study a broader class of subgroups of $\operatorname{Diff}(M)$.
The material in this section can be found in \cite{KolarMichorSlovak1993} and \cite{Saunders1989}.

\begin{defn}
	Let $M$ and $N$ be smooth manifolds, and let $r\geq 0$ be an integer.
	Two smooth curves
	$$
	\phi,\psi:\mathbb R\longrightarrow N
	$$
	are said to have \emph{$r$-th order contact at $0$}, written
	$$
	\phi\sim_r \psi,
	$$
	if, for every $f\in C^\infty(N,\mathbb R)$ and every
	$k=0, 1, \ldots, r$,
	$$
	\left.
	\frac{d^k}{dt^k}
	\right|_{t=0}
	\bigl(f\circ \phi-f\circ \psi\bigr)(t)
	=0.
	$$
	Obviously,  $\sim_r$ is an equivalence relation.
	Two smooth maps 
	$$\phi,\psi:M\longrightarrow N$$
	are said to determine the same \emph{$r$-jet at $x\in M$}
	if, 
	for every smooth curve
	$h:\mathbb R\to M$  with 
	$h(0)=x$,
	the curves  
	$\phi\circ h$ and $\psi\circ h$ have $r$-th order contact at
	zero, i.e. 
	$$
	\phi\circ h\sim_r \psi\circ h.
	$$
	In this case, we write
	$$
	j_x^r\phi=j_x^r\psi.
	$$
	The  equivalence class $j_x^r\phi$ is called the
	\emph{$r$-jet of $\phi$ at $x$}, or an \emph{$r$-jet of $M$ into $N$}.
	The set of all such $r$-jets is denoted by
	$$
	J^r(M,N).
	$$
	The points $x$ and $\phi(x)$ are called the  \emph{source} and \emph{target} of  $j_x^r\phi\in J^r(M,N)$  respectively. 
	By $J_x^r(M,N)$ or $J^r(M,N)_y$, 
	we mean the set of all
	$r$-jets of $M$ into $N$ with source $x\in M$ or target $y\in N$
	respectively, and we write
	$$
	J_x^r(M,N)_y
	=
	J_x^r(M,N)\cap J^r(M,N)_y.
	$$	
\end{defn}

The space
$$
T^{r*}M:=J^r(M,\mathbb R)_0
$$
is called the space of all \emph{$(1,r)$-covelocities} on $M$. 

Since $\mathbb R$ is a vector space,
$T^{r*}M\xrightarrow{\pi} M$
is a vector bundle, 
where $\pi(j_x^r\phi)=x$ and  each fiber $T_x^{r*}M=J_{x}^r(M,\mathbb R)_0$ 
is equipped with the operations: 
$$
j_x^r\phi+j_x^r\psi
=
j_x^r(\phi+\psi)
\quad \text{and} \quad 
a\cdot j_x^r\phi
=
j_x^r(a\cdot \phi),
\quad a\in\mathbb R.
$$
It is called  \emph{$(1,r)$-covelocities bundle on $M$}.

Consider the dual vector bundle
$$
T^{(r)}M
:=
\bigl(T^{r*}M\bigr)^*
$$
of the $(1,r)$-covelocities bundle on $M$. For every smooth map
$\phi:M\to N$, the jet composition
$$
A\longmapsto A\circ j_x^r\phi,
\quad
x\in M,\quad
A\in T_{\phi(x)}^{r*}N,
$$
defines a linear map
$$
\lambda(j_x^r\phi):
T_{\phi(x)}^{r*}N
\longrightarrow
T_x^{r*}M,
$$
and hence the dual map
$$
\bigl(\lambda(j_x^r\phi)\bigr)^*
:
T_x^{(r)}M
\longrightarrow
T_{\phi(x)}^{(r)}N.
$$
We denote
$$
D_x^{(r)}\phi = \bigl(\lambda(j_x^r\phi)\bigr)^*.
$$
They determine a smooth bundle map 
$$
D^{(r)}\phi: 
T^{(r)}M
\longrightarrow
T^{(r)}N.
$$
The assignments $M\mapsto T^{(r)}M$ and
$\phi\mapsto D^{(r)}\phi$ define a covariant functor. In particular, 
$$
D^{(r)}(\psi\circ \phi)
=
D^{(r)}\psi
\circ
D^{(r)}\phi
$$
for all smooth maps  $\phi:M_1\to M_2$  and $\psi:M_2\to M_3$.
The bundle 
$T^{(r)}M$  is called \emph{$r$-th order tangent
	bundle on $M$},
the elements of  $T^{(r)}M$ are called \emph{$r$-th order tangent
	vectors on $M$},	
$D^{(r)}\phi$
is called 
\emph{$r$-th order derivative
	map} associated to $\phi:M \to N$,
and  $D_x^{(r)}\phi$ is called the
\emph{$r$-th order derivative} of $\phi$ at $x$.

The above construction generalizes the usual tangent bundle and derivative map  since the following diagram commutes:
$$
\begin{tikzcd}[column sep=large, row sep=large]
	T^{(1)}M
	\arrow[r, "D^{(1)}\phi"]
	\arrow[d, "\cong"']
	&
	T^{(1)}N
	\arrow[d, "\cong"]
	\\
	TM
	\arrow[r, "D\phi"']
	&
	TN,
\end{tikzcd}
$$
where the vertical arrows are the canonical isomorphisms and
$D\phi:TM\to TN$ is the usual derivative map.

\subsection{Higher-order $\mu$-discreteness}

Let $T^{(r)}M$ be the $r$-th order tangent bundle of $M$. In this subsection, we first consider the metric bundle 
$$ 
\operatorname{Riem}^{(r)}(M)\longrightarrow M, 
$$ 
whose fiber over each $x\in M$ is the space of all inner products on $T_x^{(r)}M$, equipped with its canonical structure as a homogeneous Hadamard manifold. We then use the $r$-th order derivatives of diffeomorphisms to quantify how continuous sections of this bundle, viewed as Riemannian metrics on $T^{(r)}M$, are distorted. This yields a quantitative notion of higher-order metric distortion and leads naturally to the definition of \emph{higher-order $\mu$-discreteness}  for subgroups of $\operatorname{Diff}(M)$.  This notion is a generalization of that introduced by Gong--Wu--Xie--Yu \cite{GongWuXieYu}.

Let $V$ be an $n$-dimensional real vector space. Denote by 
$$ 
\mathcal{P}(V) 
:= 
\left\{ 
\textsl{g}:V\times V\longrightarrow\mathbb R ~\big|~  
\textsl{g} \text{ is a positive-definite symmetric bilinear form} 
\right\} 
$$ 
the space of all inner products on $V$. 
After choosing a linear isomorphism $V\cong\mathbb R^n$, the space 
$\mathcal{P}(V)$ is identified with 
$$ 
\mathcal{P}(n) 
:= 
\left\{ 
Q\in M_n(\mathbb R) ~\big|~ 
Q^{\mathsf T}=Q 
\ \text{and}\ 
v^{\mathsf T}Qv>0 
\text{ for every }v\neq0 
\right\}. 
$$ 
The congruence action of $\operatorname{GL}(n,\mathbb R)$ on $\mathcal{P}(n)$ 
is transitive, and the stabilizer of the identity matrix $I_n$ is 
$\operatorname O(n)$. Consequently, 
$$ 
\mathcal{P}(V) 
\cong 
\mathcal{P}(n) 
\cong 
\operatorname{GL}(n,\mathbb R)\big/ 
\operatorname O(n). 
$$ 
The standard Riemannian metric on $\mathcal{P}(n)$ is given by 
$$ 
\langle A,B\rangle_Q 
= 
\operatorname{Tr} 
\bigl(Q^{-1}AQ^{-1}B\bigr), 
\qquad 
A,B\in T_Q\mathcal{P}(n) 
\cong 
\operatorname{Sym}_n(\mathbb R). 
$$ 
This metric induces a canonical metric $d_{\mathcal{P}(V)}$ on $\mathcal{P}(V)$, 
independent of the chosen identification $V\cong\mathbb R^n$. 
With this metric, $\mathcal{P}(V)$ is a complete simply connected Riemannian 
manifold of nonpositive curvature,  
and hence a homogeneous Hadamard manifold.

Let $L:V\to V'$ be a linear isomorphism. For each
$\textsl{g}\in\mathcal P(V')$, its pullback by $L$ is the
bilinear form $L^*\textsl{g}$ on $V$ defined by
$$
(L^*\textsl{g})(u,v)
:=
\textsl{g}(Lu,Lv),
\qquad u,v\in V.
$$
This form is symmetric and positive-definite: for every
$u\in V\setminus\{0\}$, we have $Lu\neq0$ and hence
$$
(L^*\textsl{g})(u,u)
=
\textsl{g}(Lu,Lu)>0.
$$
Thus pullback defines a bijection
$$
L^*:\mathcal P(V')\longrightarrow\mathcal P(V),
\qquad
\textsl{g}\mapsto  L^*\textsl{g},
$$
whose inverse is $(L^{-1})^*$.
By the congruence invariance of the canonical metrics on
the spaces of inner products, this bijection is an isometry.
Explicitly, for all $\textsl{g},\textsl{g}'\in\mathcal P(V')$,
$$
d_{\mathcal P(V)}
\bigl(L^*\textsl{g},L^*\textsl{g}'\bigr)
=
d_{\mathcal P(V')}
\bigl(\textsl{g},\textsl{g}'\bigr).
$$

Let $M$ be a smooth manifold of dimension $n$,  and let 
$r\geq 1$ be an integer.  
The  rank of the  $r$-th order tangent bundle 
$$ 
T^{(r)}M 
= 
\bigl(J^r(M,\mathbb R)_0\bigr)^* 
$$ 
is 
$$ 
m_r 
:= 
\operatorname{rank}\bigl(T^{(r)}M\bigr) 
= 
\sum_{k=1}^{r} 
\binom{n+k-1}{k} 
= 
\binom{n+r}{r}-1. 
$$ 
By a smooth local frame of $T^{(r)}M$ over an open subset
$U\subseteq M$, we mean a smooth vector-bundle isomorphism
$$
\varepsilon:U\times\mathbb R^{m_r}
\longrightarrow
T^{(r)}M\big|_U
$$
covering $\operatorname{id}_U$. For each $x\in U$, write
$$
\varepsilon_x:\mathbb R^{m_r}\longrightarrow T_x^{(r)}M,
\qquad
\varepsilon_x(v):=\varepsilon(x,v).
$$
For $\textsl{g}\in\mathcal P(T_x^{(r)}M)$, its pullback by
$\varepsilon_x$ is the inner product
$\varepsilon_x^*\textsl{g}\in\mathcal P(\mathbb R^{m_r})$
given by
$$
(\varepsilon_x^*\textsl{g})(u,v)
=
\textsl{g}(\varepsilon_xu,\varepsilon_xv),
\qquad
u,v\in\mathbb R^{m_r}.
$$
These pullbacks define local trivializations of the metric bundle
of Hadamard manifolds
$$
\operatorname{Riem}^{(r)}(M)
:=
\bigsqcup_{x\in M}\mathcal P\bigl(T_x^{(r)}M\bigr)
\longrightarrow M,
$$
whose fiber over $x\in M$ is
$$
\operatorname{Riem}_x^{(r)}(M)
=
\mathcal P\bigl(T_x^{(r)}M\bigr).
$$
A section
$$
s\in\prod_{x\in M}\operatorname{Riem}_x^{(r)}(M)
$$
is continuous if and only if, for every smooth local frame
$\varepsilon$ over an open subset $U\subseteq M$, the map
$$
U\longrightarrow\mathcal P(\mathbb R^{m_r}),
\qquad
x\mapsto\varepsilon_x^*\bigl(s(x)\bigr)
$$
is continuous.

Now, we consider how a continuous section of this bundle, viewed as a Riemannian metric on $T^{(r)}M$, is distorted under the $r$-th order derivatives of diffeomorphisms.

Fix a continuous section  
$\textsl{g}_r: M\to \operatorname{Riem}^{(r)}(M)$ and suppose that $M$ is compact with a finite Borel measure $\mu$. 
For a linear map 
$$ 
A:T_x^{(r)}M\longrightarrow T_y^{(r)}M, 
$$ 
its operator norm with respect to  $\textsl{g}_r$ is defined by 
$$ 
\|A\|_{\textsl{g}_r} 
:= 
\sup_{u\in T_x^{(r)}M\setminus\{0\}} 
\frac{ 
	\sqrt{\textsl{g}_{r,y}(Au,Au)} 
}{ 
	\sqrt{\textsl{g}_{r,x}(u,u)} 
}. 
$$ 
For $\phi,\psi\in\operatorname{Diff}(M)$ and $x\in M$, define 
$$ 
\kappa_{r,\textsl{g}_r}(\phi,\psi,x) 
:= 
\max 
\left\{ 
\left\| 
D_{\phi^{-1}(x)}^{(r)} 
\bigl(\psi^{-1}\circ \phi\bigr) 
\right\|_{\textsl{g}_r}, 
\left\| 
D_{\psi^{-1}(x)}^{(r)} 
\bigl(\phi^{-1}\circ \psi\bigr) 
\right\|_{\textsl{g}_r} 
\right\}. 
$$ 
The two linear maps appearing here are mutually inverse, since 
$$ 
D_{\psi^{-1}(x)}^{(r)} 
\bigl(\phi^{-1}\circ \psi\bigr) 
\circ 
D_{\phi^{-1}(x)}^{(r)} 
\bigl(\psi^{-1}\circ \phi\bigr) 
= 
\operatorname{id}_{T_{\phi^{-1}(x)}^{(r)}M}. 
$$ 
Consequently, 
$$ 
\kappa_{r,\textsl{g}_r}(\phi,\psi,x)\geq 1. 
$$

Fix $1\leq p<\infty$ throughout the following discussion.

Define
$$
d_{\mu,\textsl{g}_r}^{(r,p)}:
\operatorname{Diff}(M)\times\operatorname{Diff}(M)
\longrightarrow
[0,\infty)
$$
by
$$
d_{\mu,\textsl{g}_r}^{(r,p)}(\phi,\psi)
:=
\left(
\int_M
\bigl(
\log\kappa_{r,\textsl{g}_r}(\phi,\psi,x)
\bigr)^p
\,d\mu(x)
\right)^{1/p}.
$$
The integral is well-defined.
In fact, the two maps
$$
x\longmapsto
D_{\phi^{-1}(x)}^{(r)}\bigl(\psi^{-1}\circ \phi\bigr)
\quad\text{and}\quad
x\longmapsto
D_{\psi^{-1}(x)}^{(r)}\bigl(\phi^{-1}\circ \psi\bigr)
$$
depend continuously on $x$. Since $\textsl{g}_r$ is continuous and the
operator norm of a linear map between finite-dimensional inner-product
spaces depends continuously on both the linear map and the inner
products, the function
$$
x\longmapsto\kappa_{r,\textsl{g}_r}(\phi,\psi,x)
$$
is continuous. Hence, since $M$ is compact and $\mu$ is finite,
it follows that
$$
\log\kappa_{r,\textsl{g}_r}(\phi,\psi,\cdot)\in L^p(M,\mu).
$$
Therefore, $d_{\mu,\textsl{g}_r}^{(r,p)}(\phi,\psi)$ is finite and well-defined.

\begin{lem}
	$d_{\mu,\textsl{g}_r}^{(r,p)}$ is a pseudometric on $\operatorname{Diff}(M)$,
	i.e., for any $\phi_1,\phi_2,\phi_3\in\operatorname{Diff}(M)$, we have
	\begin{itemize}
		\item [(1)] $d_{\mu,\textsl{g}_r}^{(r,p)}(\phi_1,\phi_1)=0$;
		\item [(2)] $d_{\mu,\textsl{g}_r}^{(r,p)}(\phi_1,\phi_2)
		=
		d_{\mu,\textsl{g}_r}^{(r,p)}(\phi_2,\phi_1)$;
		\item [(3)] $d_{\mu,\textsl{g}_r}^{(r,p)}(\phi_1,\phi_3)
		\leq
		d_{\mu,\textsl{g}_r}^{(r,p)}(\phi_1,\phi_2)
		+
		d_{\mu,\textsl{g}_r}^{(r,p)}(\phi_2,\phi_3)$.
	\end{itemize}
\end{lem}
\begin{proof}
	The first two axioms are immediate from the definition of
	$d_{\mu,\textsl{g}_r}^{(r,p)}$.
	To verify the last axiom, we observe that
	for
	$\phi_1,\phi_2,\phi_3\in\operatorname{Diff}(M)$, functoriality gives
	$$
	\begin{aligned}
		D_{\phi_1^{-1}(x)}^{(r)}
		\bigl(\phi_3^{-1}\circ \phi_1\bigr)
		&=
		D_{\phi_2^{-1}(x)}^{(r)}
		\bigl(\phi_3^{-1}\circ \phi_2\bigr)
		\circ
		D_{\phi_1^{-1}(x)}^{(r)}
		\bigl(\phi_2^{-1}\circ \phi_1\bigr),
		\\
		D_{\phi_3^{-1}(x)}^{(r)}
		\bigl(\phi_1^{-1}\circ \phi_3\bigr)
		&=
		D_{\phi_2^{-1}(x)}^{(r)}
		\bigl(\phi_1^{-1}\circ \phi_2\bigr)
		\circ
		D_{\phi_3^{-1}(x)}^{(r)}
		\bigl(\phi_2^{-1}\circ \phi_3\bigr).
	\end{aligned}
	$$
	The submultiplicativity of operator norms therefore implies
	$$
	\kappa_{r,\textsl{g}_r}(\phi_1,\phi_3,x)
	\leq
	\kappa_{r,\textsl{g}_r}(\phi_1,\phi_2,x)\cdot
	\kappa_{r,\textsl{g}_r}(\phi_2,\phi_3,x).
	$$
	Hence
	$$
	\big|\log\kappa_{r,\textsl{g}_r}(\phi_1,\phi_3,x)\big|
	\leq
	\log\kappa_{r,\textsl{g}_r}(\phi_1,\phi_2,x)
	+
	\log\kappa_{r,\textsl{g}_r}(\phi_2,\phi_3,x),
	$$
	and the Minkowski inequality yields
	$$
	d_{\mu,\textsl{g}_r}^{(r,p)}(\phi_1,\phi_3)
	\leq
	d_{\mu,\textsl{g}_r}^{(r,p)}(\phi_1,\phi_2)
	+
	d_{\mu,\textsl{g}_r}^{(r,p)}(\phi_2,\phi_3).
	$$
\end{proof}

For $\tau\in\operatorname{Diff}(M)$, the pseudometric satisfies the
covariance relation
$$
d_{\tau_*\mu,\textsl{g}_r}^{(r,p)}(\phi,\psi)
=
d_{\mu,\textsl{g}_r}^{(r,p)}
\bigl(\tau^{-1}\circ \phi,\tau^{-1}\circ \psi\bigr),
$$
where $\tau_*\mu$ denotes the push-forward of $\mu$ under $\tau$,
defined by
$$
(\tau_*\mu)(E):=\mu\bigl(\tau^{-1}(E)\bigr)
$$
for every Borel set $E\subseteq M$.

We now restrict the pseudometric to a countable subgroup
$\Gamma<\operatorname{Diff}(M)$. For every $\gamma\in\Gamma$,
$$
\begin{aligned}
	d_{\mu,\textsl{g}_r}^{(r,p)}
	\bigl(\gamma,1_{\Gamma}\bigr)
	=
	\Bigg(
	\int_M
	\Bigg(
	\log\max\Bigg\{
	&
	\left\|
	D_{\gamma^{-1}(x)}^{(r)}\gamma
	\right\|_{\textsl{g}_r},
	\\
	&
	\left\|
	D_x^{(r)}(\gamma^{-1})
	\right\|_{\textsl{g}_r}
	\Bigg\}
	\Bigg)^p
	\,d\mu(x)
	\Bigg)^{1/p}.
\end{aligned}
$$
Since
$$
D_x^{(r)}(\gamma^{-1})
=
\left(
D_{\gamma^{-1}(x)}^{(r)}\gamma
\right)^{-1},
$$
this pseudometric measures, in an $L^p$-sense, how far the $r$-th order
derivative of $\gamma$ is from being an isometry.

We introduce the following higher-order generalization of
$\mu$-discreteness given in \cite{GongWuXieYu}.

\begin{defn}
	We say that a countable subgroup
	$$
	\Gamma<\operatorname{Diff}(M)
	$$
	is \emph{$(r,p)$-th order $\mu$-discrete} if
	$$
	\inf_{\tau\in\Gamma}
	d_{\tau_*\mu,\textsl{g}_r}^{(r,p)}
	\bigl(\gamma,1_\Gamma\bigr)
	\longrightarrow
	\infty
	\qquad
	\text{as }\gamma\to\infty\text{ in }\Gamma.
	$$
	Equivalently, for every $R>0$, the set
	$$
	\left\{
	\gamma\in\Gamma:
	\inf_{\tau\in\Gamma}
	d_{\tau_*\mu,\textsl{g}_r}^{(r,p)}(\gamma,1_\Gamma)
	\leq R
	\right\}
	$$
	is finite.
\end{defn}

\begin{rmk}
	The definition of  $(r,p)$-th  order $\mu$-discreteness is independent
	of the choice of the continuous section $\textsl{g}_r$.
	Indeed, let $\textsl{g}'_r$ be another continuous section.
	Since $M$ is compact, there exists $C\geq1$ such that
	$$
	C^{-1}\textsl{g}_r
	\leq
	\textsl{g}'_r
	\leq
	C\textsl{g}_r
	$$
	as quadratic forms on each fiber. Consequently,
	$$
	C^{-1}\kappa_{r,\textsl{g}_r}(\phi,\psi,x)
	\leq
	\kappa_{r,\textsl{g}'_r}(\phi,\psi,x)
	\leq
	C\kappa_{r,\textsl{g}_r}(\phi,\psi,x)
	$$
	for all $\phi,\psi\in\operatorname{Diff}(M)$ and $x\in M$.
	Taking logarithms and applying the reverse triangle inequality
	in $L^p$ gives
	$$
	\left|
	d_{\nu,\textsl{g}_r}^{(r,p)}(\phi,\psi)
	-
	d_{\nu,\textsl{g}'_r}^{(r,p)}(\phi,\psi)
	\right|
	\leq
	\nu(M)^{1/p}\log C
	$$
	for every finite Borel measure $\nu$ on $M$.
	Since $(\tau_*\mu)(M)=\mu(M)$, this bound is uniform for
	$\nu=\tau_*\mu$, $\tau\in\Gamma$.
	Thus the divergence condition in the definition is unchanged
	when $\textsl{g}_r$ is replaced by $\textsl{g}'_r$.
	
	Moreover, it is straightforward to verify that, when $r=1$ and $p=2$,
	our definition reduces to the notion of $\mu$-discreteness
	introduced in \cite{GongWuXieYu}.
\end{rmk}

\subsection{Higher-order orbit maps}

Throughout this subsection, we fix $1\leq p<\infty$ and let
$M$ be a compact smooth manifold equipped with a finite
Borel measure $\mu$.
We construct the higher-order orbit map from $\operatorname{Diff}(M)$
into the $L^p$-direct integral of the metric bundle
$\operatorname{Riem}^{(r)}(M)$.
We then prove that, for any countable subgroup
$\Gamma<\operatorname{Diff}(M)$ equipped with a proper
left-invariant metric, the restriction of this orbit map to
$\Gamma$ is a coarse embedding if and only if
$\Gamma$ is $(r,p)$-th order $\mu$-discrete.
At last, we prove Corollary \ref{maintheoremofkdjalkf}.

\begin{constr}
	Assume  that $M$ is compact and that $\mu$ is a finite
	Borel measure on $M$. We fix a base
	section 
	$$
	\textsl{g}_r:M\longrightarrow\operatorname{Riem}^{(r)}(M), \quad x \mapsto \textsl{g}_{r,x}
	$$
	that is continuous.
	Let  $\phi\in\operatorname{Diff}(M)$. 
	Since
	$$
	D^{(r)}(\phi^{-1}\circ \phi)
	=
	D^{(r)}\operatorname{id}_M=\operatorname{id}_{T^{(r)}M},
	$$
	the linear map
	$$
	D_x^{(r)}\phi:
	T_x^{(r)}M
	\longrightarrow
	T_{\phi(x)}^{(r)}M
	$$
	is an isomorphism whose inverse is
	$$
	D_{\phi(x)}^{(r)}(\phi^{-1}).
	$$
	For every section
	$$
	s:M\longrightarrow\operatorname{Riem}^{(r)}(M),
	\quad
	x\mapsto s_x,
	$$
	of $\operatorname{Riem}^{(r)}(M)$, define its
	\emph{$r$-th order push-forward} by $\phi$ to be the section
	$\phi_*^{(r)}s$ given by
	$$
	\bigl(\phi_*^{(r)}s\bigr)_y
	:=
	\bigl(D_y^{(r)}(\phi^{-1})\bigr)^*
	s_{\phi^{-1}(y)}.
	$$
	Equivalently, for $u,v\in T_y^{(r)}M$,
	$$
	\bigl(\phi_*^{(r)}s\bigr)_y(u,v)
	=
	s_{\phi^{-1}(y)}
	\left(
	D_y^{(r)}(\phi^{-1})(u),
	D_y^{(r)}(\phi^{-1})(v)
	\right).
	$$
	In particular, if $y=\phi(x)$, then
	$$
	\bigl(\phi_*^{(r)}s\bigr)_{\phi(x)}
	=
	\left(
	\bigl(D_x^{(r)}\phi\bigr)^{-1}
	\right)^*s_x.
	$$
	Since pullback by a linear isomorphism preserves positive-definiteness,
	$\bigl(\phi_*^{(r)}s\bigr)_y$ belongs to
	$\mathcal P(T_y^{(r)}M)$. Hence $\phi_*^{(r)}s$ is again a section of
	$\operatorname{Riem}^{(r)}(M)$.
	
	The higher-order push-forward satisfies
	$$
	\operatorname{id}_*^{(r)}s=s
	$$
	and
	$$
	(\phi\circ \psi)_*^{(r)}s
	=
	\phi_*^{(r)}\bigl(\psi_*^{(r)}s\bigr)
	$$
	for all $\phi,\psi\in\operatorname{Diff}(M)$. 
	
	Moreover, this push-forward preserves continuity and
	Borel measurability.
	Indeed, let $\varepsilon$ be a smooth local frame over $U\subseteq M$.
	It induces a smooth local frame over $\phi(U)$ by
	$$
	\widetilde{\varepsilon}_y
	:=
	D_{\phi^{-1}(y)}^{(r)}\phi
	\circ
	\varepsilon_{\phi^{-1}(y)},
	\qquad y\in\phi(U).
	$$
	By the definition of $\phi_*^{(r)}$, 
	the local representation of $\phi_*^{(r)}s$ in this frame is
	$$
	\begin{aligned}
		\phi(U)&\longrightarrow\mathcal P(\mathbb R^{m_r}),\\
		y&\longmapsto
		\widetilde{\varepsilon}_y^*
		\bigl((\phi_*^{(r)}s)_y\bigr)
		=
		\varepsilon_{\phi^{-1}(y)}^*
		\bigl(s_{\phi^{-1}(y)}\bigr).
	\end{aligned}
	$$
	This is the composition of the local representation
	$x\mapsto\varepsilon_x^*(s_x)$ with $\phi^{-1}$.
	Since $\phi^{-1}$ is continuous, $\phi_*^{(r)}s$ is continuous
	whenever $s$ is continuous, and Borel measurable whenever
	$s$ is Borel measurable.

	Since the pointwise distance between continuous sections is continuous,
	compactness of $M$ and finiteness of $\mu$ give
	$$
	\int_M
	d_{\mathcal P(T_x^{(r)}M)}
	\left(
	\bigl(\phi_*^{(r)}\textsl{g}_r\bigr)_x,
	\textsl{g}_{r,x}
	\right)^p
	\,d\mu(x)
	<\infty.
	$$
	Consequently,
	$$
	\bigl[\phi_*^{(r)}\textsl{g}_r\bigr]
	\in
	L^p_{\textsl{g}_r}
	\bigl(M,\mu,\operatorname{Riem}^{(r)}(M)\bigr)
	$$
	for every $\phi\in\operatorname{Diff}(M)$.
	We define the \emph{$(r,p)$-th  order orbit map}
	$$
	F_{\mu,\textsl{g}_r}^{(r,p)}:
	\operatorname{Diff}(M)
	\longrightarrow
	L^p_{\textsl{g}_r}
	\bigl(M,\mu,\operatorname{Riem}^{(r)}(M)\bigr)
	$$
	by
	$$
	F_{\mu,\textsl{g}_r}^{(r,p)}(\phi)
	:=
	\bigl[\phi_*^{(r)}\textsl{g}_r\bigr].
	$$
\end{constr}

We are now ready to give an equivalent characterization of higher-order $\mu$-discreteness in terms of the higher-order orbit map.

\begin{prop}\label{whatiskfdsjlkjsa}
	Assume that $M$ is compact smooth manifold and $\mu$ is a finite
	Borel measure on $M$. 
	For any countable subgroup $\Gamma<\operatorname{Diff}(M)$ equipped with a proper left-invariant metric $d_{\Gamma}$,	
	the following statements are equivalent:
	\begin{itemize}
		\item [(a)]  $\Gamma$ is  $(r,p)$-th  order $\mu$-discrete;
		\item [(b)] 	$F_{\mu,\textsl{g}_r}^{(r,p)}\big|_\Gamma$ is a coarse embedding for every continuous section 
		$$\textsl{g}_r:M\longrightarrow\operatorname{Riem}^{(r)}(M).$$
	\end{itemize}
\end{prop}
\begin{proof}
	The finite case is immediate, so assume that $\Gamma$ is infinite.
	
	Fix a continuous base section $\textsl{g}_r$ and write
	$$
	F:=F_{\mu,\textsl{g}_r}^{(r,p)}\big|_\Gamma.
	$$
	For $\gamma,\widetilde\gamma\in\Gamma$, set
	$$
	\psi:=\widetilde\gamma^{-1}\gamma,
	\qquad
	\tau:=\widetilde\gamma^{-1}.
	$$
	By functoriality,
	$$
	\tau_*^{(r)}
	\bigl(\gamma_*^{(r)}\textsl{g}_r\bigr)
	=
	\psi_*^{(r)}\textsl{g}_r,
	\qquad
	\tau_*^{(r)}
	\bigl(\widetilde\gamma_*^{(r)}\textsl{g}_r\bigr)
	=
	\textsl{g}_r.
	$$
	Hence, we obtain
	\begin{equation}\label{eq:diyigedongxikldjiikjsdo}
    \begin{aligned}	&d_{L^p}\bigl(F(\gamma),F(\widetilde\gamma)\bigr)^p
		\\
		&\quad=
		\int_M
		d_{\mathcal P(T_y^{(r)}M)}
		\left(
		\bigl(\gamma_*^{(r)}\textsl{g}_r\bigr)_y,
		\bigl(\widetilde\gamma_*^{(r)}\textsl{g}_r\bigr)_y
		\right)^p
		\,d\mu(y)
		\\
		&\quad=
		\int_M
		d_{\mathcal P(T_{\tau(y)}^{(r)}M)}
		\left(
		\bigl(\tau_*^{(r)}
		(\gamma_*^{(r)}\textsl{g}_r)\bigr)_{\tau(y)},
		\bigl(\tau_*^{(r)}
		(\widetilde\gamma_*^{(r)}\textsl{g}_r)\bigr)_{\tau(y)}
		\right)^p
		\,d\mu(y)
		\\
		&\quad=
		\int_M
		d_{\mathcal P(T_{\tau(y)}^{(r)}M)}
		\left(
		\bigl(\psi_*^{(r)}\textsl{g}_r\bigr)_{\tau(y)},
		\textsl{g}_{r,\tau(y)}
		\right)^p
		\,d\mu(y)
		\\
		&\quad=
		\int_M
		d_{\mathcal P(T_x^{(r)}M)}
		\left(
		\bigl(\psi_*^{(r)}\textsl{g}_r\bigr)_x,
		\textsl{g}_{r,x}
		\right)^p
		\,d(\tau_*\mu)(x).
	\end{aligned}
\end{equation}
	We next express the integrand in terms of singular values.
	Put
	$$
	m_r:=\operatorname{rank}\bigl(T^{(r)}M\bigr).
	$$
	For each $x\in M$, choose a
	$\textsl{g}_{r,\psi^{-1}(x)}$-orthonormal basis of
	$T_{\psi^{-1}(x)}^{(r)}M$ and a
	$\textsl{g}_{r,x}$-orthonormal basis of $T_x^{(r)}M$.
	Let $A_x$ be the matrix of
	$D_{\psi^{-1}(x)}^{(r)}\psi$ under these bases, and let
	$$
	\sigma_1(x)\geq\cdots\geq\sigma_{m_r}(x)>0
	$$
	be its singular values.
	These values are independent of the chosen orthonormal bases;
	no smooth choice of bases is needed.
	Since
	$$
	D_x^{(r)}(\psi^{-1})
	=
	\bigl(D_{\psi^{-1}(x)}^{(r)}\psi\bigr)^{-1},
	$$
	the definition of $\kappa_{r,\textsl{g}_r}$ gives
	\begin{equation}\label{eq:diergedongxiklsdjfljdsl}
	\begin{aligned}
		\kappa_{r,\textsl{g}_r}(\psi,1_\Gamma,x)
		&=
		\max\left\{
		\sigma_1(x),\sigma_{m_r}(x)^{-1}
		\right\},
		\\
		\log\kappa_{r,\textsl{g}_r}(\psi,1_\Gamma,x)
		&=
		\max_{1\leq j\leq m_r}
		\bigl|\log\sigma_j(x)\bigr|.
	\end{aligned}
	\end{equation}
	In the chosen basis of $T_x^{(r)}M$, the matrix of
	$\textsl{g}_{r,x}$ is $I_{m_r}$.
	By the definition of the higher-order push-forward, the matrix of
	$\bigl(\psi_*^{(r)}\textsl{g}_r\bigr)_x$ is
	$$
	Q_x
	:=
	A_x^{-\mathsf T}A_x^{-1}
	=
	\bigl(A_xA_x^{\mathsf T}\bigr)^{-1}.
	$$
	The eigenvalues of $A_xA_x^{\mathsf T}$ are
	$\sigma_1(x)^2,\ldots,\sigma_{m_r}(x)^2$, so those of $Q_x$ are
	$\sigma_1(x)^{-2},\ldots,\sigma_{m_r}(x)^{-2}$.
	The distance formula for the trace metric on $\mathcal P(m_r)$
	therefore yields
	\begin{equation}\label{eq:disangedongxidkljlk}
	\begin{aligned}
		&
		d_{\mathcal P(T_x^{(r)}M)}
		\left(
		\bigl(\psi_*^{(r)}\textsl{g}_r\bigr)_x,
		\textsl{g}_{r,x}
		\right)
		\\
		&\quad=
		d_{\mathcal P(m_r)}(Q_x,I_{m_r})
		\\
		&\quad=
		\left(
		\sum_{j=1}^{m_r}
		\bigl(\log(\sigma_j(x)^{-2})\bigr)^2
		\right)^{1/2}
		\\
		&\quad=
		2\left(
		\sum_{j=1}^{m_r}
		\bigl(\log\sigma_j(x)\bigr)^2
		\right)^{1/2}.
	\end{aligned}
	\end{equation}
	Combining \eqref{eq:diergedongxiklsdjfljdsl}  and \eqref{eq:disangedongxidkljlk},
	we obtain
	\begin{equation}\label{eq:disigedongxikjljklk}
	\begin{aligned}
		2\log\kappa_{r,\textsl{g}_r}(\psi,1_\Gamma,x)
		&\leq
		d_{\mathcal P(T_x^{(r)}M)}
		\left(
		\bigl(\psi_*^{(r)}\textsl{g}_r\bigr)_x,
		\textsl{g}_{r,x}
		\right)
		\\
		&\leq
		2\sqrt{m_r}\,
		\log\kappa_{r,\textsl{g}_r}(\psi,1_\Gamma,x).
	\end{aligned} 
	\end{equation}
	Since $\kappa_{r,\textsl{g}_r}(\psi,1_\Gamma,x)\geq1$,
	all terms in this pointwise comparison are nonnegative.
	Hence, it follows from \eqref{eq:diyigedongxikldjiikjsdo}  and \eqref{eq:disigedongxikjljklk}  that 
	$$
	\begin{aligned}
		&
		2^p d_{\tau_*\mu,\textsl{g}_r}^{(r,p)}(\psi,1_\Gamma)^p
		\\
		&\quad=
		2^p\int_M
		\bigl(
		\log\kappa_{r,\textsl{g}_r}(\psi,1_\Gamma,x)
		\bigr)^p
		\,d(\tau_*\mu)(x)
		\\
		&\quad\leq
		\int_M
		d_{\mathcal P(T_x^{(r)}M)}
		\left(
		\bigl(\psi_*^{(r)}\textsl{g}_r\bigr)_x,
		\textsl{g}_{r,x}
		\right)^p
		\,d(\tau_*\mu)(x)
		\\
		&\quad=
		d_{L^p}\bigl(F(\gamma),F(\widetilde\gamma)\bigr)^p
		\\
		&\quad\leq
		(2\sqrt{m_r})^p\int_M
		\bigl(
		\log\kappa_{r,\textsl{g}_r}(\psi,1_\Gamma,x)
		\bigr)^p
		\,d(\tau_*\mu)(x)
		\\
		&\quad=
		(2\sqrt{m_r})^p\,d_{\tau_*\mu,\textsl{g}_r}^{(r,p)}(\psi,1_\Gamma)^p.
	\end{aligned}
	$$
	Taking $p$-th roots therefore yields
	\begin{equation}\label{eq:diwugedongxidd}
	\begin{aligned}
		2d_{\tau_*\mu,\textsl{g}_r}^{(r,p)}(\psi,1_\Gamma)
		&\leq
		d_{L^p}\bigl(F(\gamma),F(\widetilde\gamma)\bigr)
		\\
		&\leq
		2\sqrt{m_r}\,
		d_{\tau_*\mu,\textsl{g}_r}^{(r,p)}(\psi,1_\Gamma).
	\end{aligned}
	\end{equation}
	
	$(b)\Longrightarrow(a)$.  
	Suppose that $F$ is a coarse embedding with lower control function
	$\rho_-$, where $\rho_-(t)\to\infty$ as $t\to\infty$.
	By \eqref{eq:diwugedongxidd}  and the left invariance of $d_{\Gamma}$,
	$$
	\begin{aligned}
		\rho_-\bigl(d_\Gamma(\psi,1_\Gamma)\bigr)
		&\leq
		d_{L^p}\bigl(F(\widetilde\gamma\psi),F(\widetilde\gamma)\bigr)\\
		&\leq
		2\sqrt{m_r}\,
		d_{\tau_*\mu,\textsl{g}_r}^{(r,p)}(\psi,1_\Gamma).
	\end{aligned}
	$$
	Thus,
	$$
	\inf_{\tau\in\Gamma}
	d_{\tau_*\mu,\textsl{g}_r}^{(r,p)}(\psi,1_\Gamma)
	\geq
	\frac{\rho_-\bigl(d_\Gamma(\psi,1_\Gamma)\bigr)}
	{2\sqrt{m_r}}
	\longrightarrow\infty
	\qquad\text{as }~ d_\Gamma(\psi,1_\Gamma)\to\infty,
	$$
	i.e., $\Gamma$ is $(r,p)$-th  order $\mu$-discrete.

	$(a)\Longrightarrow(b)$.  
	Suppose that $\Gamma$ is  $(r,p)$-th  order $\mu$-discrete.
	For $R\geq0$, put
	$$
	C_R
	:=
	\max_{\substack{\psi\in\Gamma\\ d_\Gamma(\psi,1_\Gamma)\leq R}}
	\sup_{x\in M}
	\log\kappa_{r,\textsl{g}_r}(\psi,1_\Gamma,x).
	$$
	Properness of $d_\Gamma$, compactness of $M$, and continuity of
	$\kappa_{r,\textsl{g}_r}$ imply that $C_R$ is finite and
	nondecreasing. Since $(\tau_*\mu)(M)=\mu(M)$, \eqref{eq:diwugedongxidd}  gives
	$$
	d_{L^p}\bigl(F(\gamma),F(\widetilde\gamma)\bigr)
	\leq 2\sqrt{m_r}\,\mu(M)^{1/p}\,C_R
	\quad\text{whenever }d_\Gamma(\gamma,\widetilde\gamma)\leq R.
	$$
	Thus $\rho_+(R):= 2\sqrt{m_r}\,\mu(M)^{1/p}\,C_R$ is an upper
	control function.
	
	For the lower control, define
	$$
	a(\psi)
	:=
	\inf_{\tau\in\Gamma}
	d_{\tau_*\mu,\textsl{g}_r}^{(r,p)}(\psi,1_\Gamma),
	\qquad
	\rho_-(R)
	:=
	2\inf_{\substack{\psi\in\Gamma\\ d_\Gamma(\psi,1_\Gamma)\geq R}}
	a(\psi).
	$$
	Since $\Gamma$ is infinite and $d_\Gamma$ is proper, $\rho_-$
	is finite-valued and nondecreasing. Moreover, $(r,p)$-th  order $\mu$-discreteness
	implies that $\rho_-(R)\to\infty$ as $R\to\infty$.
	For $\psi=\widetilde\gamma^{-1}\gamma$ and
	$\tau=\widetilde\gamma^{-1}$,  \eqref{eq:diwugedongxidd}  yields
	$$
	\begin{aligned}
		\rho_-\bigl(d_\Gamma(\gamma,\widetilde\gamma)\bigr)
		&\leq 2a(\psi)\\
		&\leq 2d_{\tau_*\mu,\textsl{g}_r}^{(r,p)}(\psi,1_\Gamma)\\
		&\leq d_{L^p}\bigl(F(\gamma),F(\widetilde\gamma)\bigr).
	\end{aligned}
	$$
	Hence $F$ is a coarse embedding.
\end{proof}

\section{Snowflake metric bundle and embedding into $L^p$-spaces}

In this section, we introduce the notion of a \emph{snowflake metric bundle}, whose fibers admit snowflake embeddings into a separable Hilbert space.
By the work of Lang, Schlichenmaier, Naor and Silberman \cite{LangSchlichenmaier,NaorSilberman2011}, we know that this class includes every metric bundle whose typical fiber is a separable metric space of finite Nagata dimension. Examples of such fibers include homogeneous Hadamard manifolds, pinched negatively curved Hadamard manifolds, doubling spaces, separable metric trees, and separable Euclidean buildings, as well as finite products and finite unions of these spaces. In particular, the higher-order Riemannian metric bundles discussed in Section~3 provide natural examples of snowflake metric bundles.
We then show that the associated $L^p$-direct integrals (spaces of measurable sections of snowflake metric bundles) admit coarse embeddings into $L^{p'}$-spaces for suitable $p'\in[1,\infty)$. Consequently, the Novikov conjecture (respectively, the coarse Baum--Connes conjecture) holds for every countable group (respectively, discrete metric space with bounded geometry) that admits a coarse embedding into one of these direct integrals.

\subsection{$L^p$-direct integral and snowflake metric bundle}

In this subsection, we introduce metric bundles and their $L^p$-direct integrals. We then focus on snowflake metric bundles.

Let 
$$
\pi:E\longrightarrow\Omega.
$$ 
be a continuous surjective map between topological spaces.
Assume that each fiber $E_z:=\pi^{-1}(z)$ carries a metric $d_z$
compatible with its subspace topology.
We call $\pi:E\to\Omega$ a \emph{metric bundle}
if there exists a metric space $(X,d_X)$ such that every point
of $\Omega$ has an open neighborhood $U$ with a homeomorphism
$$
\Phi_U:\pi^{-1}(U)\longrightarrow U\times X,
$$
whose restriction to $(E_z,d_z)$ is an isometry onto
$\{z\}\times X$ for every $z\in U$.
We call $\Phi_U$ a \emph{local trivialization}
and $(X,d_X)$ the \emph{typical fiber}.

A metric bundle $\pi:E\to\Omega$ is called \emph{measurable}
if its  \emph{fiberwise distance function}
$$
\delta:E\times_\Omega E\longrightarrow[0,\infty),
\qquad
\delta(u,v):=d_{\pi(u)}(u,v),
$$
is measurable with respect to the restriction of
$\mathcal B(E)\otimes\mathcal B(E)$ to
$$
E\times_\Omega E
:=
\{(u,v)\in E\times E:\pi(u)=\pi(v)\},
$$
where $\mathcal B(E)$ denotes the Borel $\sigma$-algebra of $E$, and
$\mathcal B(E)\otimes\mathcal B(E)$ is the product $\sigma$-algebra.

Let  $\mu$  be a 
Borel measure  
on $\Omega$. 
A section $s:\Omega\to E$  is called
\emph{$\mu$-measurable} if $s^{-1}(U)$ is $\mu$-measurable
for every open subset $U\subseteq E$.
The measurability of the fiberwise distance function implies that,
for any two $\mu$-measurable sections $s,t$, the function
$$
z\longmapsto d_z\bigl(s(z),t(z)\bigr)
$$
is $\mu$-measurable.

\begin{defn}\label{generaldedifofproduct} 
	Let 
	$\pi:E\to\Omega$ 
	be a measurable metric bundle,  and  let  $\mu$  be a
	Borel measure  
	on $\Omega$. 
	Fix a $\mu$-measurable base section $s_0:\Omega\to E$.
	For $p\in[1,+\infty)$,
	the \emph{$L^p$-direct integral} of $E$ over $(\Omega,\mu)$
	based at $s_0$, denoted by
	$$
	L^p_{s_0}(\Omega,\mu,E),
	$$
	is the space of equivalence classes of $\mu$-measurable
	sections $s$ satisfying
	$$
	\int_\Omega
	d_z\bigl(s(z),s_0(z)\bigr)^p\,d\mu(z)<\infty.
	$$
	Two sections are equivalent if and only if they agree
	$\mu$-almost everywhere.
	Moreover, the Minkowski inequality
	implies  that  the formula 
	$$
	d_{L^p}([s],[t])
	=
	\left(
	\int_\Omega
	d_z\bigl(s(z),t(z)\bigr)^p\,d\mu(z)
	\right)^{1/p}
	$$
	defines a well-defined metric on $L^p_{s_0}(\Omega,\mu,E)$.
\end{defn}

\begin{rmk}
	We  consider two  important special cases.
	\begin{itemize}
		\item [(1)]  Let $\Omega=M$  be a smooth manifold equipped with
		a $\sigma$-finite    Borel measure $\mu$, and set
		$E=\operatorname{Riem}^{(r)}(M)$,
		the metric bundle over $M$
		introduced in Section~3.2.
		Roughly speaking, its fiber over each $x\in M$ consists of all
		inner products on the $r$-th order tangent space $T_x^{(r)}M$
		and carries a canonical Riemannian metric that makes it a
		homogeneous Hadamard manifold.

		With the canonical fiber metrics
		$d_{\mathcal P(T_x^{(r)}M)}$, $x\in M$,
		the fiberwise distance function is continuous.
		Thus, for any $\mu$-measurable base section $s_0$, the preceding
		definition gives the $L^p$-direct integral
		$$L^p_{s_0}\bigl(M,\mu, \operatorname{Riem}^{(r)}(M)\bigr).
		$$
		Note that the $L^2$-direct integral is a Hilbert--Hadamard space in the framework of Gong-Wu-Xie-Yu \cite{GongWuXieYu}.  We avoid using this concept for the sake of brevity. 
		
		\item [(2)]  
		For a trivial metric bundle $\Omega\times X\to\Omega$,
		the preceding definition gives a based version of the
		$L^p$-continuum product
		(for $p=2$ with finite measure and a constant base section,
		see \cite[Construction~2.19]{GongWuXieYu}).
		If, in addition, $X$ is a separable Banach space
		equipped with its norm metric and the base section
		is the zero section, the preceding definition recovers
		the classical Bochner space.
		
		More precisely, let $\Omega$ be a topological space equipped
		with a Borel measure $\mu$, let $(X,d_X)$ be a separable
		metric space, and let $p\in[1,\infty)$.
		Consider the trivial bundle $\Omega\times X\to\Omega$,
		whose fibers are naturally identified with $(X,d_X)$.
		Since $X$ is separable, the fiberwise distance function is
		automatically measurable.
		Fix a $\mu$-measurable base section
		$s_0:\Omega\to\Omega\times X$, and write
		$$
		s_0(z)=(z,f_0(z)),
		$$
		where $f_0:\Omega\to X$ is $\mu$-measurable.
		Under the identification of $\mu$-measurable sections
		$s(z)=(z,f(z))$ with $\mu$-measurable maps $f:\Omega\to X$,
		the space
		$$
		L^p_{s_0}(\Omega,\mu,\Omega\times X)
		$$
		coincides isometrically with $L^p_{f_0}(\Omega,\mu,X)$.
		Here $L^p_{f_0}(\Omega,\mu,X)$  is called  the
		\emph{$L^p$-continuum product} of $X$ over $(\Omega,\mu)$
		based at $f_0$.
		It consists of equivalence classes  of
		$\mu$-measurable maps $f:\Omega\to X$ satisfying
		$$
		\int_\Omega d_X\bigl(f(z),f_0(z)\bigr)^p\,d\mu(z)<\infty,
		$$
		and is equipped with the metric
		$$
		d_{L^p}([f],[g])
		=
		\left(
		\int_\Omega d_X\bigl(f(z),g(z)\bigr)^p\,d\mu(z)
		\right)^{1/p}.
		$$
		If, in addition, $X$ is a Banach space equipped with
		the norm metric $d_X(x,y)=\|x-y\|$, then its separability
		ensures that every $\mu$-measurable map $\Omega\to X$
		is Bochner measurable
		(see \cite[Appendix~K]{AlbiacKalton2016}).
		Thus, for the zero base section $s_0(z)=(z,0)$,
		the space $L^p_{s_0}(\Omega,\mu,\Omega\times X)$
		(i.e., $L^p_0(\Omega,\mu,X)$) is canonically isometric to
		the classical Bochner space $L^p(\Omega,\mu,X)$.
	\end{itemize}
\end{rmk}

We next introduce the notion of snowflake embedding.
Let $(X,d_X)$ and $(Y,d_Y)$ be metric spaces.
For $\alpha\in(0,1]$, a map $h:X\to Y$ is called an
\emph{$\alpha$-snowflake embedding} if
$$
h:(X,d_X^\alpha)\longrightarrow(Y,d_Y)
$$
is a bi-Lipschitz embedding.
Equivalently, there exist constants $0<c\leq C<\infty$
such that
$$
c\,d_X(x,y)^\alpha
\leq
d_Y\bigl(h(x),h(y)\bigr)
\leq
C\,d_X(x,y)^\alpha
\quad \text{for all }  x,y\in X.
$$
If such a map exists, we say that $X$ admits an
\emph{$\alpha$-snowflake embedding}  into $Y$.

The study of snowflake embeddings has roots in Schoenberg's
work in the 1930s \cite{Schoenberg1938}.
A major subsequent development was Assouad's classical theorem,
which established that every doubling metric space admits
an $\alpha$-snowflake embedding into a finite-dimensional
Euclidean space for every $\alpha\in(0,1)$
\cite{Assouad1983}.

The theory later expanded beyond the doubling setting.
Lang and Schlichenmaier studied Nagata dimension, a quantitative
refinement of Gromov's asymptotic dimension, and showed that
metric spaces of finite Nagata dimension admit bi-Lipschitz
embeddings into finite products of metric trees after snowflaking
by sufficiently small exponents
\cite[Theorem~1.3]{LangSchlichenmaier}.
Naor and Silberman subsequently showed, through a more subtle
argument, that every metric space of finite Nagata dimension
admits an $\alpha$-snowflake embedding into Hilbert space
for every $\alpha\in(0,1)$
\cite[Lemma~5.1 and Theorem~5.2]{NaorSilberman2011}.

The class of metric spaces with finite Nagata dimension includes  all homogeneous Hadamard manifolds,
pinched negatively curved Hadamard manifolds,
doubling spaces,  geodesic Gromov hyperbolic spaces that are
doubling in the small,  metric ($\mathbb{R}$-)trees, and  Euclidean
buildings, among others. 
Moreover, this class is closed under taking finite
products and finite unions \cite{LangSchlichenmaier}.

These results motivate our introduction of snowflake metric bundles.
By applying snowflake embeddings fiberwise and then assembling
them, we obtain coarse embeddings of the associated
$L^p$-direct integrals into suitable Banach spaces, leading to
applications to the coarse Baum--Connes and Novikov conjectures.
Bundles arising from higher-order Riemannian metrics, which we
discuss later, provide natural examples  of this framework.

\begin{defn}
	Let $\pi:E\to\Omega$ be a 
	metric bundle with   typical fiber $(X,d_X)$. 
	The bundle is call a  \emph{snowflake metric bundle}
	if  $X$ admits an $\alpha$-snowflake embedding
	into separable Hilbert space for some  $\alpha \in (0,1]$.
\end{defn}

\begin{rmk}
The class of snowflake metric bundles includes metric bundles whose
typical fibers have finite Nagata dimension and therefore
covers a broad range of examples from metric geometry.
The associated $L^p$-direct integrals include classical
$L^p$-spaces and a broad class of nonlinear analogues,
especially those arising from the higher-order Riemannian
metric bundles introduced later.
\end{rmk}

\subsection{Embedding into Banach spaces}

In this subsection, we show  that  $L^p$-direct integrals of snowflake metric bundles
over metacompact spaces equipped with $\sigma$-finite Borel
measures admit coarse embeddings into $L^{p'}$-spaces for
suitable $p'\in[1,\infty)$.

\begin{lem}\label{prop:metric-bundle-trivialization}
	Let $\Omega$ be a metacompact space equipped with a Borel
	measure $\mu$, and let $\pi:E\to\Omega$ be a 
	metric bundle with a separable typical fiber $(X,d_X)$. 	Then  $E$
	is a  measurable metric bundle.
	Moreover, for every  $p\in[1,+\infty)$ and  $\mu$-measurable base section $s_0:\Omega\to E$,
	there exists  an isometric bijection
	$$
	\widehat{\Psi}:
	L^p_{s_0}(\Omega,\mu,E)
	\longrightarrow
	L^p_{\widehat{s_0}}(\Omega,\mu, X),
	$$
	where  $\widehat{s_0}:\Omega\to  X$  is a $\mu$-measurable map associated to  $s_0$.
\end{lem}
\begin{proof}
	By metacompactness, there exists a point-finite open cover
	$\{U_i\}_{i\in\Lambda}$ of $\Omega$
	(that is, each point belongs to only finitely many $U_i$)
	with local trivializations
	$$
	\Phi_i:\pi^{-1}(U_i)\longrightarrow U_i\times X.
	$$
	For each $z\in\Omega$, let
	$$
	\Lambda(z):=\{i\in\Lambda:z\in U_i\}.
	$$
	This set is finite and nonempty.
	
	For $n\geq1$, put
	$$
	G_n:=\{z\in\Omega:|\Lambda(z)|\geq n\},
	\qquad
	B_n:=G_n\setminus G_{n+1}.
	$$
	Every set $G_n$ is open, since
	$$
	G_n=
	\bigcup_{\substack{F\subseteq\Lambda\\|F|=n}}
	\bigcap_{i\in F}U_i.
	$$
	Thus each $B_n$ is Borel.
	Since every $\Lambda(z)$ is finite and nonempty,
	$$
	\Omega=\bigsqcup_{n=1}^{\infty}B_n.
	$$
	Hence the  sets $B_n$ form a countable Borel
	partition of $\Omega$, even though $\Lambda$ may be
	uncountable.
	
	By Zermelo's well-ordering theorem, fix a well-order on
	$\Lambda$.
	For each $z\in\Omega$, let $i(z)$ be the least element
	of $\Lambda(z)$ with respect to this order.
	Define
	$$
	\Psi:E\longrightarrow\Omega\times X,
	\qquad
	\Psi(u):=\Phi_{i(\pi(u))}(u).
	$$
	This is a bijection preserving base points, with inverse
	$$
	\Psi^{-1}(z,x):=\Phi_{i(z)}^{-1}(z,x).
	$$
	Let $\widehat{\pi}:\Omega\times X\to X$ denote the projection
	onto the second factor.
	By construction,
	\begin{equation}\label{eq:innerproductofhat}
	d_X\bigl(\widehat{\pi}(\Psi(u)),\widehat{\pi}(\Psi(v))\bigr)
	=
	d_z(u,v),
	\qquad  u,v\in E_z.
	\end{equation}
	
	We claim that, for every $n\geq1$, the restriction
	$$
	\Psi_n:=\Psi|_{\pi^{-1}(B_n)}:
	\pi^{-1}(B_n)\longrightarrow B_n\times X
	$$
	is a homeomorphism.
	Indeed, for $z\in B_n$, the set
	$$
	V_z:=\bigcap_{i\in\Lambda(z)}U_i
	$$
	is an open neighborhood of $z$.
	If $z'\in V_z\cap B_n$, then
	$\Lambda(z)\subseteq\Lambda(z')$ and both sets have
	exactly $n$ elements.
	Thus $\Lambda(z')=\Lambda(z)$, so $i(z')=i(z)$.
	Consequently,
	$$
	\Psi|_{\pi^{-1}(V_z\cap B_n)}:
	\pi^{-1}(V_z\cap B_n)
	\longrightarrow
	(V_z\cap B_n)\times X
	$$
	is the restriction of $\Phi_{i(z)}$ and is therefore
	a homeomorphism.
	As $z$ ranges over $B_n$, the sets
	$\pi^{-1}(V_z\cap B_n)$ and $(V_z\cap B_n)\times X$
	form open covers of $\pi^{-1}(B_n)$ and $B_n\times X$,
	respectively, in their subspace topologies.
	Hence both $\Psi_n$ and $\Psi_n^{-1}$ are continuous,
	proving the claim.
	
	The sets $\pi^{-1}(B_n)$ and $B_n\times X$ are Borel,
	since $\pi$ and the projection $\Omega\times X\to\Omega$
	are continuous.
	They give countable Borel decompositions of $E$ and
	$\Omega\times X$, respectively.
	As each $\Psi_n$ is a homeomorphism, both $\Psi$
	and $\Psi^{-1}$ are Borel measurable.
	
	Since $X$ is separable, the function
	$(x,y)\mapsto d_X(x,y)$ is measurable with respect to
	the product $\sigma$-algebra
	$\mathcal B(X)\otimes\mathcal B(X)$, and hence the fiberwise
	distance function of $\Omega\times X$ is measurable.
	Together with the Borel measurability of $\Psi$ and its
	preservation of distances within each fiber, this implies
	that the fiberwise distance function of $E$ is also measurable.
	
	For a  section $s:\Omega\to E$, set
	$$
	\widehat{s}:=\widehat{\pi}\circ\Psi\circ s.
	$$
	Since $X$ is separable, a map $f:\Omega\to X$ is
	$\mu$-measurable if and only if the associated section
	$z\mapsto(z,f(z))$ of the trivial bundle
	$\Omega\times X\to\Omega$ is $\mu$-measurable.
	Together with the Borel measurability of $\Psi$ and
	$\Psi^{-1}$, this shows that
	$$
	s \longmapsto \widehat{s}
	$$
	is a bijection between $\mu$-measurable sections of $E$
	and $\mu$-measurable maps $\Omega\to X$.
	Moreover,
	$$
	s=t\quad\mu\text{-a.e.}
	\quad\Longleftrightarrow\quad
	\widehat{s}=\widehat{t}\quad\mu\text{-a.e.}
	$$
	Combining  \eqref{eq:innerproductofhat}  with  these observations, we obtain
	$$
	\int_\Omega
	d_z\bigl(s(z),t(z)\bigr)^p\,d\mu(z)
	=
	\int_\Omega
	d_X\bigl(\widehat{s}(z),\widehat{t}(z)\bigr)^p\,d\mu(z)
	$$
	for any two $\mu$-measurable sections $s,t:\Omega\to E$.
	It therefore induces a well-defined isometric bijection
	$$
	\widehat{\Psi}:
	L^p_{s_0}(\Omega,\mu,E)
	\longrightarrow
	L^p_{\widehat{s_0}}(\Omega,\mu,X),
	\qquad
	\widehat{\Psi}([s]):=[\widehat{s}].
	$$
\end{proof}

We emphasize that every snowflake metric bundle is measurable
by Lemma \ref{prop:metric-bundle-trivialization},
since its typical fiber is separable.
Consequently, the construction of $L^p$-direct integrals
in Definition \ref{generaldedifofproduct}
applies to every snowflake metric bundle.
The following result is a generalization of Proposition 12.3 in \cite{TianYu2026}.

\begin{prop}\label{snowflakemetricbundleembedding}
	Let $\Omega$ be a metacompact space equipped with a $\sigma$-finite Borel
	measure $\mu$, and let $\pi:E\to\Omega$ be a 
	snowflake metric bundle  with   typical fiber $(X,d_X)$. 
	Then for every $p\in[1,\infty)$,
	there exists $p'\in[1,\infty)$ such that the
	$L^p$-direct integral
	$$
	L^p_{s_0}(\Omega,\mu,E)
	$$
	admits a coarse embedding into an $L^{p'}$-space
	$L^{p'}(\Omega',\mu')$   for any  $\mu$-measurable base section $s_0:\Omega\to E$.
\end{prop}
\begin{proof}
	Fix $p\in[1,\infty)$.
	
	Since $X$ admits a snowflake embedding into a separable
	Hilbert space, $X$ is separable.
	By Lemma  \ref{prop:metric-bundle-trivialization}, it therefore
	suffices to show that
	$$
	L^p_{f_0}(\Omega,\mu,X)
	$$
	admits a coarse embedding into  an $L^{p'}$-space
	$L^{p'}(\Omega',\mu')$,
	where $f_0: \Omega\to  X$  is  a $\mu$-measurable map.  
	
	Choose a snowflake embedding
	$$
	h:X\longrightarrow\ell^2(\mathbb N).
	$$
	Thus there exist $\alpha\in(0,1]$ and constants
	$0<c\leq C<\infty$ such that
	\begin{equation}\label{eq:inequalityofhwith}
	c\,d_X(x,y)^\alpha
	\leq
	\|h(x)-h(y)\|
	\leq
	C\,d_X(x,y)^\alpha,
	\qquad  x,y\in X.
	\end{equation}
	Set
	$$
	p':=\frac{p}{\alpha}.
	$$
	Then $p'\in[1,\infty)$ and $\alpha p'=p$.
	
	Let $m$ denote Lebesgue measure on $[0,1]$, and put
	$$
	Y:=L^{p'}([0,1],m).
	$$
	The Gaussian construction  (see Proposition 6.4.12 in \cite{AlbiacKalton2016})  gives a linear isometric embedding
	$$
	\mathcal{I}:\ell^2(\mathbb N)\longrightarrow Y.
	$$
	Indeed, let $\{g_n\}_{n\geq1}$ be independent standard real Gaussian
	random variables on $([0,1],m)$.
	These variables belong to $L^{p'}([0,1],m)$, 
	and every finite real linear combination
	$\sum_j a_jg_j$ has the same distribution as
	$(\sum_j a_j^2)^{1/2}g_1$.
	Hence
	$$
	\left\|\sum_j a_jg_j\right\|_{L^{p'}}
	=
	\big\|g_1\big\|_{L^{p'}}
	\left(\sum_j |a_j|^2\right)^{1/2}.
	$$
	Let  $\{e_n\}_{n\geq1}$  denote the canonical orthonormal basis
	of $\ell^2(\mathbb N)$.
	The assignment
	$e_n\mapsto g_n/\|g_1\|_{L^{p'}}$
	extends  to a linear isometric embedding
	$\mathcal{I}: \ell^2(\mathbb{N}) \to  Y$.
	
	For a $\mu$-measurable map $f:\Omega\to X$ representing
	an element of $L^p_{f_0}(\Omega,\mu,X)$, define
	$\widetilde{f}:\Omega\to Y$ by
	$$
	\widetilde{f}(z)
	:=
	\mathcal I\Bigl(h\bigl(f(z)\bigr)-h\bigl(f_0(z)\bigr)\Bigr),
	\qquad z\in\Omega.  
	$$
	The continuity of $h$ and $\mathcal I$, together with the
	separability of $Y$, implies that $\widetilde{f}$
	is Bochner measurable with respect to $\mu$.
	Moreover, it follows from \eqref{eq:inequalityofhwith}  that
	$$
	\begin{aligned}
		\int_\Omega  \big\| \widetilde{f}(z) \big\|_Y^{p'}\,d\mu(z)
		&=
		\int_\Omega
		\Big\| h\big(f(z)\big)-h\big(f_0(z)\big) \Big\|_2^{p'}\,d\mu(z)
		\\
		&\leq
		C^{p'}
		\int_\Omega
		d_X\big( f(z), f_0(z) \big)^p\,d\mu(z)
		<\infty.
	\end{aligned}
	$$
	Thus $\widetilde{f}$ represents an element of the Bochner
	space $L^{p'}(\Omega,\mu, Y)$.
	Note that the class $[\widetilde{f}]$ depends only on $[f]$.
	Hence the map
	$$
	\psi:
	L^p_{f_0}(\Omega,\mu,X)
	\longrightarrow L^{p'}(\Omega,\mu, Y),
	\qquad
	\psi([f]):=[\widetilde{f}],
	$$
	is well-defined.
	
	For any $[f],[g]\in L^p_{f_0}(\Omega,\mu,X)$,
	the linearity and isometry of $\mathcal I$ give
	$$
	\begin{aligned}
		\Big\| \psi([f])-\psi([g]) \Big\|_{L^{p'}(\Omega,\mu,Y)}^{p'}
		&=
		\int_\Omega
		\big\| \widetilde{f}(z)-\widetilde{g}(z) \big\|_Y^{p'}\,d\mu(z)\\
		&=
		\int_\Omega
		\Big\|h\big(f(z)\big)-h\big(g(z)\big)\Big\|_2^{p'}
		\,d\mu(z).
	\end{aligned}
	$$
	It follows from  \eqref{eq:inequalityofhwith}  that
	$$
	c\,d_{L^p}([f],[g])^\alpha
	\leq
	\Big\| \psi([f])-\psi([g]) \Big\|_{L^{p'}(\Omega,\mu,Y)}
	\leq
	C\,d_{L^p}([f],[g])^\alpha.
	$$
	Thus, $\psi$ is a coarse embedding. 
	
	Finally, set
	$$
	\Omega':=\Omega\times[0,1] \quad \text{and} 
	\quad
	\mu':=\mu\otimes m,
	$$
	where $\mu\otimes m$  denotes the product measure.
	By the Fubini isomorphism for Bochner spaces
	(see \cite[Proposition~1.2.24]{HNVW2016}), there is a canonical isometric isomorphism
	$$
	L^{p'}(\Omega,\mu, Y)
	\xlongrightarrow{\cong} L^{p'}(\Omega',\mu').
	$$
	Composing $\psi$ with this isomorphism  gives the
	required coarse embedding.
	This completes the proof.
\end{proof}

\subsection{Proof of the main theorem}

We first briefly recall the definitions of  the strong Novikov conjecture and the coarse Baum--Connes conjecture for  countable discrete group. 
 
Let $\Gamma$ be a countable discrete group equipped with a proper left-invariant metric, and denote the so-obtained metric space as $|\Gamma|$.

For $R>0$, let $P_R(\Gamma)$ denote the
Rips complex with vertex set $\Gamma$, in which a finite subset spans
a simplex if and only if its diameter is at most $R$. The Rips complex $P_R(\Gamma)$ is equipped with the aspherical metric (cf.~ \cite{Y2000}).

The group $\Gamma$ acts on $P_R(\Gamma)$ by $\gamma\cdot(\sum_{g}a_{g}g)=\sum_{g}a_{g}\gamma g$ for $\sum_{g}a_{g}\gamma g$ and $\gamma\in \Gamma$. This is a proper action by isometries.

The Baum--Connes and coarse Baum--Connes assembly maps are,
respectively,
\begin{align*}
\mu_\Gamma \colon
\varinjlim_{R\to\infty} K_*^\Gamma(P_R(\Gamma))
&\longrightarrow K_*(C_r^*(\Gamma)),\\
\mu_X \colon
\varinjlim_{R\to\infty} K_*(P_R(X))
&\longrightarrow K_*(C^*(X)).
\end{align*}
Here $\Gamma$ acts on $P_R(\Gamma)$ by left translation,
$K_*^\Gamma(P_R(\Gamma))$ denotes its equivariant $K$-homology,
and $K_*(P_R(X))$ denotes the locally finite $K$-homology
of $P_R(X)$.
Both direct limits are induced by inclusions of Rips complexes
as $R$ increases.
Moreover, $C_r^*(\Gamma)$ is the reduced group $C^*$-algebra
of $\Gamma$, and $C^*(X)$ is the Roe algebra of $X$.
The \emph{strong Novikov conjecture} for $\Gamma$ and the
\emph{coarse Baum--Connes conjecture} for $X$ assert,
respectively, that $\mu_\Gamma$ is injective and $\mu_X$ is isomorphism
for $*=0,1$.
Kasparov and Yu \cite{KasparovYu2012} proved that the Novikov conjecture holds for every countable discrete group admitting a coarse embedding into an $\ell^p$-space for some $1\leq p<\infty$.
Wang,  Xie, Yu and  Zhu \cite{wangxieyuzhu}  showed that the coarse Baum-Connes conjecture holds for any discrete metric space with boundary geometry that admits a coarse embedding into an $\ell^p$-space for some $1\leq p<\infty$.

We need the following useful lemma.

\begin{lem}\label{prop:Lp-to-ellp}
	Let $1\leq p<\infty$, and let $(\Gamma,d_\Gamma)$ be either
	a countable discrete group equipped with a proper left-invariant
	metric or a discrete metric space with bounded geometry.
	If $(\Gamma,d_\Gamma)$ admits a coarse embedding into an
	$L^p$-space, then it admits a coarse embedding into
	$\ell^p(\mathbb N)$.
\end{lem}
\begin{proof}

	Suppose that
	$$
	h:\Gamma\longrightarrow L^p(\Omega,\mu)
	$$
	is a coarse embedding. Thus, there exist nondecreasing functions
	$$
	\rho_-,\rho_+:[0,\infty)\longrightarrow[0,\infty)
	$$
	such that
	$\lim_{t\to\infty}\rho_-(t)=\infty$
	and
	$$
	\rho_-\bigl(d_\Gamma(\gamma,\widetilde{\gamma})\bigr)
	\leq
	\|h(\gamma)-h(\widetilde{\gamma})\|_{L^p}
	\leq
	\rho_+\bigl(d_\Gamma(\gamma,\widetilde{\gamma})\bigr)
	$$
	for all $\gamma,\widetilde{\gamma}\in\Gamma$.
	
	We first show that all finite subsets of $\Gamma$ admit uniformly
	controlled coarse embeddings into $\ell^p(\mathbb N)$. Let
	$E\subseteq\Gamma$ be finite. Since simple functions are dense in
	$L^p(\Omega,\mu)$, for every $\gamma\in E$ we may choose a simple
	function $s_\gamma\in L^p(\Omega,\mu)$ such that
	$$
	\|h(\gamma)-s_\gamma\|_{L^p}<\frac12.
	$$
	Because $E$ is finite, after passing to a common refinement of the
	finitely many measurable partitions involved, there exist pairwise
	disjoint measurable sets
	$$
	A_1,\ldots,A_N\subseteq\Omega,
	\qquad
	0<\mu(A_j)<\infty,
	$$
	and scalars $c_{\gamma,j}$ such that
	$$
	s_\gamma
	=
	\sum_{j=1}^{N}c_{\gamma,j}\mathbf 1_{A_j}
	\qquad
	\text{for every }\gamma\in E.
	$$
	Define
	$$
	F_E:E\longrightarrow\ell_N^p\subseteq\ell^p(\mathbb N)
	$$
	by
	$$
	F_E(\gamma)
	:=
	\bigl(
	\mu(A_1)^{1/p}c_{\gamma,1},
	\ldots,
	\mu(A_N)^{1/p}c_{\gamma,N}
	\bigr).
	$$
	For all $\gamma,\widetilde{\gamma}\in E$, we then have
	$$
	\begin{aligned}
		\|F_E(\gamma)-F_E(\widetilde{\gamma})\|_{\ell^p}^p
		&=
		\sum_{j=1}^{N}
		\mu(A_j)
		\lvert
		c_{\gamma,j}-c_{\widetilde{\gamma},j}
		\rvert^p \\
		&=
		\|s_\gamma-s_{\widetilde{\gamma}}\|_{L^p}^p.
	\end{aligned}
	$$
	Hence
	$$
	\|F_E(\gamma)-F_E(\widetilde{\gamma})\|_{\ell^p}
	=
	\|s_\gamma-s_{\widetilde{\gamma}}\|_{L^p}.
	$$
	Moreover, the reverse triangle inequality gives
	$$
	\begin{aligned}
		\Bigl|
		\|s_\gamma-s_{\widetilde{\gamma}}\|_{L^p}
		-
		\|h(\gamma)-h(\widetilde{\gamma})\|_{L^p}
		\Bigr|
		&\leq
		\bigl\|
		(s_\gamma-h(\gamma))
		-
		(s_{\widetilde{\gamma}}-h(\widetilde{\gamma}))
		\bigr\|_{L^p} \\
		&\leq
		\|s_\gamma-h(\gamma)\|_{L^p}
		+
		\|s_{\widetilde{\gamma}}-h(\widetilde{\gamma})\|_{L^p} \\
		&<1.
	\end{aligned}
	$$
	Consequently,
	$$
	\max\left\{
	0,\,
	\rho_-\bigl(
	d_\Gamma(\gamma,\widetilde{\gamma})
	\bigr)-1
	\right\}
	\leq
	\|F_E(\gamma)-F_E(\widetilde{\gamma})\|_{\ell^p}
	\leq
	\rho_+\bigl(
	d_\Gamma(\gamma,\widetilde{\gamma})
	\bigr)+1.
	$$
	Define
	$$
	\widetilde{\rho}_-(t)
	:=
	\max\{0,\rho_-(t)-1\},
	\qquad
	\widetilde{\rho}_+(t)
	:=
	\rho_+(t)+1.
	$$
	These functions are independent of the finite subset $E$, and
	$$
	\lim_{t\to\infty}\widetilde{\rho}_-(t)=\infty.
	$$
	Thus, all finite subsets of $\Gamma$ admit uniformly controlled
	coarse embeddings into $\ell^p(\mathbb N)$.
	
In either case, every bounded subset of $\Gamma$ is finite.
We may apply the finite-determination theorem
for coarse embeddings of locally finite metric spaces into
Banach spaces
\cite[Theorem~1.3]{OstrovskiiFiniteDetermination}
to conclude that $\Gamma$ admits a coarse embedding into
$\ell^p(\mathbb N)$.
\end{proof}

We now present a technical theorem, which will lead to the proof of the main theorem.

\begin{thm}\label{maintheorem1}
	Let $\Omega$ be a metacompact space equipped with a $\sigma$-finite Borel
	measure $\mu$, and let $\pi:E\to\Omega$ be a 
	snowflake metric bundle.
Then  the Novikov conjecture
(respectively, the coarse Baum--Connes conjecture) holds for every countable  group
(respectively, discrete metric space with bounded geometry)
that admits a coarse embedding into
the
$L^p$-direct integral
$$
L^p_{s_0}(\Omega,\mu,E)
$$
for some $\mu$-measurable base section $s_0:\Omega\to E$   and some $p\in[1,+\infty)$.
\end{thm}
\begin{proof}
Let $\Gamma$ be a countable discrete group equipped with a proper
left-invariant metric (respectively, a discrete metric space
with bounded geometry). Suppose that $\Gamma$ admits a coarse
embedding into the $L^p$-direct integral
$$
L^p_{s_0}(\Omega,\mu,E)
$$
for some $p\in[1,\infty)$ and some $\mu$-measurable base section
$s_0:\Omega\to E$.

Combining Proposition \ref{snowflakemetricbundleembedding}
and Lemma \ref{prop:Lp-to-ellp}, we obtain a coarse embedding of
$\Gamma$ into $\ell^{p'}(\mathbb{N})$ for some
$p'\in[1,\infty)$.

The desired conclusion now follows from the results of
Kasparov--Yu \cite{KasparovYu2012} and
Wang--Xie--Yu--Zhu \cite{wangxieyuzhu}.
\end{proof}

Finally, we are ready to prove the Theorem \ref{main-theorem}.

\begin{proof}[Proof of Theorem \ref{main-theorem}]
	
The desired result follows immediately from Proposition \ref{whatiskfdsjlkjsa}  and	Theorem \ref{maintheorem1}.  
\end{proof}

\begin{rmk}
	Note that the coarse embeddability of the group $\Gamma$ into
	$$
	L^p_{\textsl{g}_r}
	\bigl(M,\mu,\operatorname{Riem}^{(r)}(M)\bigr)
	$$
	is a weaker condition than $(r,p)$-th order $\mu$-discreteness,
	since the latter requires the specific orbit map
	$F_{\mu,\textsl{g}_r}^{(r,p)}\big|_\Gamma$
	to be a coarse embedding.
	Consequently, the corollary applies to a broader class of groups
	than those covered by the higher-order $\mu$-discreteness criterion.
\end{rmk}

\section{New examples}

Throughout this section, we fix $p=2$ and write
$F_{\mu,\textsl{g}_r}^{(r)}$ for
$F_{\mu,\textsl{g}_r}^{(r,2)}$.
We present two families of examples demonstrating that higher-order
$\mu$-discreteness has strictly greater detecting power than its
first-order counterpart. 
We first consider amenable subgroups of $\operatorname{Diff}(M)$
and then turn to non-amenable subgroups.
In both families, the first-order $\mu$-discreteness method does
not apply, whereas our higher-order method establishes the Novikov
conjecture for the groups under consideration.

By Proposition \ref{whatiskfdsjlkjsa}, $\Gamma$ is $r$-th order
$\mu$-discrete if and only if the restricted orbit map
$F_{\mu,\textsl{g}_r}^{(r)}\big|_\Gamma$ is a coarse embedding.
We therefore formulate the  examples in terms of these orbit maps.

\begin{ex}\label{ex:heisenberg-second-order}
	There exists a finitely generated amenable subgroup
	$$
	\Gamma<\operatorname{Diff}(S^2)
	$$
	with the following properties:
	\begin{enumerate}
		\item [(a)]
		for every  finite  Borel   measure $\mu$ on $S^2$
		and every continuous section
		$$
		\textsl{g}_1:S^2\longrightarrow\operatorname{Riem}^{(1)}(S^2),
		$$
		the orbit map
		$$
		F_{\mu,\textsl{g}_1}^{(1)}\big|_\Gamma:
		\Gamma
		\longrightarrow
		L^2_{\textsl{g}_1}
		\bigl(
		S^2,\mu,\operatorname{Riem}^{(1)}(S^2)
		\bigr),
		\qquad
		F_{\mu,\textsl{g}_1}^{(1)}(\gamma)
		=
		\bigl[\gamma_*^{(1)}\textsl{g}_1\bigr],
		$$
		is not a coarse embedding;
		
		\item [(b)]
		there exists a    Borel probability measure $\mu_0$ on $S^2$ such
		that  for every  continuous section
		$$
		\textsl{g}_2:S^2\longrightarrow\operatorname{Riem}^{(2)}(S^2),
		$$
		the orbit map		
		$$
		F_{\mu_0,\textsl{g}_2}^{(2)}\big|_\Gamma:
		\Gamma
		\longrightarrow
		L^2_{\textsl{g}_2}
		\bigl(
		S^2,\mu_0,\operatorname{Riem}^{(2)}(S^2)
		\bigr),
		\qquad
		F_{\mu_0,\textsl{g}_2}^{(2)}(\gamma)
		=
		\bigl[\gamma_*^{(2)}\textsl{g}_2\bigr],
		$$
		is a coarse embedding.
	\end{enumerate}
\end{ex}
\begin{proof}
	Let
	$$
	\Gamma
	:=
	\left\{
	g_{(a,b,c)}
	=
	\begin{pmatrix}
		1&a&c\\
		0&1&b\\
		0&0&1
	\end{pmatrix}
	:
	a,b,c\in\mathbb Z
	\right\}
	$$
	be the discrete Heisenberg group. 
	Its multiplication law  is
	$$
	g_{(a,b,c)}\cdot  g_{(a',b',c')}
	=
	g_{(a+a',b+b',c+c'+ab')}.
	$$
	It is a non-abelian  nilpotent
	group generated by 
	$$
	X:=g_{(1,0,0)}
	\quad \text{and} \quad 
	Y:=g_{(0,1,0)}.
	$$
	In particular, $\Gamma$ is amenable.
	For the finite symmetric generating set
	$$
	\Delta:=\{X,X^{-1},Y,Y^{-1}\},
	$$
	we equip $\Gamma$ with the associated word metric $d_\Gamma$.

	Every matrix $A\in\operatorname{GL}(3,\mathbb R)$ determines a
	smooth map
	$$
	\phi_A:S^2\longrightarrow S^2,
	\qquad
	\phi_A(v):=\frac{Av}{\|Av\|}.
	$$
	Since
	$\phi_A\circ\phi_B=\phi_{AB}$  for all $A,B\in\operatorname{GL}(3,\mathbb R)$ and $\phi_{I}=\operatorname{id}_{S^2}$, 
	we obtain a group homomorphism
	$$
	\operatorname{GL}(3,\mathbb R)\longrightarrow\operatorname{Diff}(S^2), \quad  A\longmapsto\phi_A.
	$$
	The kernel of this homomorphism is  $\{tI: t>0\}$.
	Since
	$$
	\Gamma\cap\{tI:t>0\}=\{I\},
	$$
	the restriction of this homomorphism to $\Gamma$ is injective. 
	We henceforth identify $\Gamma$ with its image in
	$\operatorname{Diff}(S^2)$.

	{\bf  Step one}. 
	We show that	for every  finite  Borel  measure $\mu$ on $S^2$
	and every continuous section
	$\textsl{g}_1:S^2\to\operatorname{Riem}^{(1)}(S^2)$,
	the orbit  map
	$F_{\mu,\textsl{g}_1}^{(1)}$
	is not a coarse embedding.
	
	Set
	$$
	p:=(1,0,0)\in S^2.
	$$
	Every element of $\Gamma$ fixes both $p$ and $-p$. Consider the open
	hemisphere
	$$
	U_+:=\{(x_1,x_2,x_3)\in S^2:x_1>0\}
	$$
	and the coordinate parametrization
	$$
	\kappa:\mathbb R^2\longrightarrow U_+,
	\qquad
	\kappa(y,z)
	=
	\frac{(1,y,z)}{\sqrt{1+y^2+z^2}}.
	$$
	For any $g_{(a,b,c)}\in\Gamma$, the local expression of $\phi_{g_{(a,b,c)}}$  near $p$  is
	$$
	\varPhi_{a,b,c}
	:=
	\kappa^{-1}\circ\phi_{g_{(a,b,c)}}\circ\kappa,
	$$
	where
	$$
	\varPhi_{a,b,c}(y,z)
	=
	\left(
	\frac{y+bz}{1+ay+cz},
	\frac{z}{1+ay+cz}
	\right).
	$$
	Since
	$$
	D_0\varPhi_{a,b,c}
	=
	\begin{pmatrix}
		1&b\\
		0&1
	\end{pmatrix},
	$$
	it follows that  the  derivative of $\phi_{g_{(a,b,c)}}$  at $p$
	detects only the parameter $b$.
	
	In particular,
	$$
	D_0\varPhi_{a,0,0}=I,
	$$
	and hence
	$$
	D_p^{(1)}\phi_{g_{(a,0,0)}}
	=
	\operatorname{id}_{T_p^{(1)}S^2}
	$$
	for all $a\in\mathbb{Z}$.
	Similarly,
	$$
	D_{-p}^{(1)}\phi_{g_{(a,0,0)}}
	=
	\operatorname{id}_{T_{-p}^{(1)}S^2}
	$$
	for all $a\in\mathbb{Z}$.

	Fix a    Borel probability measure $\mu$ on $S^2$ and a
	continuous section
	$$
	\textsl{g}_1:S^2\longrightarrow\operatorname{Riem}^{(1)}(S^2).
	$$
	For every $a\in\mathbb Z$, define
	$$
	\delta_{a,\textsl{g}_1}(x)
	:=
	d_{\mathcal P(T_x^{(1)}S^2)}
	\bigg(
	\Bigl(\big(\phi_{g_{(a,0,0)}}\big)_*^{(1)}\textsl{g}_1\Bigr)_x,
	~ \textsl{g}_{1,x}
	\bigg).
	$$
	The function $\delta_{a,\textsl{g}_1}$ is continuous on $S^2$. Since  $\phi_{g_{(a,0,0)}}$  fixes
	$p$ and $-p$, and its first-order derivatives at these points are the
	identity, we have
	\begin{equation}
	\delta_{a,\textsl{g}_1}(p)
	=
	\delta_{a,\textsl{g}_1}(-p)
	=
	0.
	\tag{1}
	\end{equation}
	
	For each integer $m\geq1$, set
	$$
	\tau_m:=g_{(m,m,m^3)}.
	$$
	If $x=(x_1,x_2,x_3)\in S^2$, then
	$$
	\phi_{\tau_m}(x)
	=
	\dfrac{
		\left(
		x_1+mx_2+m^3x_3,\,
		x_2+mx_3,\,
		x_3
		\right)
	}{
		\sqrt{
			\left(x_1+mx_2+m^3x_3\right)^2
			+\left(x_2+mx_3\right)^2
			+x_3^2
		}
		~}.
	$$
	It follows that
	\begin{equation*}
    \phi_{\tau_m}(x)
	\longrightarrow
	\begin{cases}
		\operatorname{sgn}(x_3)p,
		& x_3\neq0,\\
		\operatorname{sgn}(x_2)p,
		&x_3=0,\ x_2\neq 0,\\
		\operatorname{sgn}(x_1)p,
		&x_2=x_3=0,
	\end{cases} \qquad  \text{as } m\to\infty.
	\end{equation*}
    
	In every case, the limit belongs to $\{p,-p\}$. 
	Combining this fact with equation $(1)$ and the continuity of $\delta_{a,\textsl{g}_1}$,
	we obtain, for each fixed $a\in\mathbb Z$ and $x\in S^2$,
	$$
	\delta_{a,\textsl{g}_1}\bigl(\phi_{\tau_m}(x)\bigr)
	\longrightarrow 0
	\qquad\text{as } m\to\infty.
	$$
	Hence, for each fixed  $a\in\mathbb{Z}$,  the dominated convergence theorem implies that
	$$
	\int_{S^2}
	\delta_{a,\textsl{g}_1}\bigl(\phi_{\tau_m}(x)\bigr)^2
	\,d\mu(x)
	\longrightarrow 0
	\qquad
	\text{as }m\to\infty.
	$$
	Now,
	for each integer $a\geq1$, choose $m_a$ sufficiently
	large such that
	\begin{equation}
	\int_{S^2}
	\delta_{a,\textsl{g}_1}
	\bigl(\phi_{\tau_{m_a}}(x)\bigr)^2
	\,d\mu(x)
	\leq \frac{1}{a^2}.
	\tag{2}
	\end{equation}
	Set
	$$
	\alpha_a:=\tau_{m_a}^{-1}g_{(a,0,0)}
	\quad\text{and}\quad
	\beta_a:=\tau_{m_a}^{-1}.
	$$
	Since $g_{(a,0,0)}=X^a$ and
	$d_\Gamma(X^a,1_\Gamma)=a$
	for every $a\geq1$, the left invariance of $d_\Gamma$ gives
	\begin{equation}
	\begin{aligned}
		d_\Gamma(\alpha_a,\beta_a)
		&=
		d_\Gamma\bigl(
		\tau_{m_a}^{-1}g_{(a,0,0)},
		\tau_{m_a}^{-1}
		\bigr)\\
		&=
		d_\Gamma\bigl(g_{(a,0,0)},1_\Gamma\bigr)\\
		&=
		d_\Gamma(X^a,1_\Gamma)
		=a
		\longrightarrow\infty
		\qquad\text{as }a\longrightarrow+\infty. 
	\end{aligned} \tag{3} 
	\end{equation}

	We next compute the distance between their images under the orbit map
	$F_{\mu,\textsl{g}_1}^{(1)}$, from which we deduce that
	$F_{\mu,\textsl{g}_1}^{(1)}\big|_\Gamma$ is not a coarse embedding.

	For every $x\in S^2$, the fiberwise isometry induced by
	$\phi_{\beta_a^{-1}}$ gives
	$$
	\begin{aligned}
		&
		d_{\mathcal P(T_x^{(1)}S^2)}
		\left(
		\bigl((\phi_{\alpha_a})_*^{(1)}\textsl{g}_1\bigr)_x,~
		\bigl((\phi_{\beta_a})_*^{(1)}\textsl{g}_1\bigr)_x
		\right)
		\\
		&\quad=
		d_{\mathcal P(T_{\phi_{\beta_a^{-1}}(x)}^{(1)}S^2)}
		\left(
		\left(
		(\phi_{\beta_a^{-1}})_*^{(1)}
		(\phi_{\alpha_a})_*^{(1)}\textsl{g}_1
		\right)_{\phi_{\beta_a^{-1}}(x)},~
		\left(
		(\phi_{\beta_a^{-1}})_*^{(1)}
		(\phi_{\beta_a})_*^{(1)}\textsl{g}_1
		\right)_{\phi_{\beta_a^{-1}}(x)}
		\right)
		\\
		&\quad=
		d_{\mathcal P(T_{\phi_{\beta_a^{-1}}(x)}^{(1)}S^2)}
		\left(
		\bigl((\phi_{g_{(a,0,0)}})_*^{(1)}\textsl{g}_1\bigr)
		_{\phi_{\beta_a^{-1}}(x)},~
		\textsl{g}_{1,\phi_{\beta_a^{-1}}(x)}
		\right)
		\\
		&\quad=
		\delta_{a,\textsl{g}_1}\bigl(\phi_{\beta_a^{-1}}(x)\bigr)
		\\
		&\quad=
		\delta_{a,\textsl{g}_1}
		\bigl(\phi_{\tau_{m_a}}(x)\bigr).
	\end{aligned}
	$$
	Combining this equality with $(2)$, we obtain that
	$$
	\begin{aligned}
		&
		d_{L^2}
		\left(
		F_{\mu,\textsl{g}_1}^{(1)}(\alpha_a),
		F_{\mu,\textsl{g}_1}^{(1)}(\beta_a)
		\right)^2
		\\
		&\quad=
		\int_{S^2}
		d_{\mathcal P(T_x^{(1)}S^2)}
		\left(
		\bigl((\phi_{\alpha_a})_*^{(1)}\textsl{g}_1\bigr)_x,~
		\bigl((\phi_{\beta_a})_*^{(1)}\textsl{g}_1\bigr)_x
		\right)^2
		\,d\mu(x)
		\\
		&\quad=
		\int_{S^2}
		\delta_{a,\textsl{g}_1}
		\bigl(\phi_{\tau_{m_a}}(x)\bigr)^2
		\,d\mu(x)
		\\
		&\quad\leq \frac{1}{a^2}
		\longrightarrow 0
		\qquad\text{as }a\longrightarrow+\infty.
	\end{aligned}
	$$
	Combining this with $(3)$, we conclude that
	$F_{\mu,\textsl{g}_1}^{(1)}\big|_\Gamma$ is not a coarse embedding.
	
	{\bf Step two}.	
	We show that 
	there exists a    Borel probability measure $\mu_0$ on $S^2$ 
	such that for every
	continuous section
	$\textsl{g}_2:S^2\to\operatorname{Riem}^{(2)}(S^2)$,
	$F_{\mu_0,\textsl{g}_2}^{(2)}\big|_\Gamma$
	is a coarse embedding.

	Since every $\gamma\in\Gamma$ fixes $p$, the second-order derivative 
	$D_p^{(2)}\phi_{\gamma}$ 
	is an endomorphism of $T_p^{(2)}S^2$. Moreover,
	functoriality of the second-order tangent map yields
	a group homomorphism
	$$
	\rho_{2,p}:
	\Gamma\longrightarrow
	\operatorname{GL}\bigl(T_p^{(2)}S^2\bigr),
	\qquad
	\rho_{2,p}(\gamma):=D_p^{(2)}\phi_{\gamma}.
	$$
	
We next show that $\rho_{2,p}$ is an injective homomorphism. This means
	that, for a given element $g_{(a,b,c)}\in\Gamma$, the second-order derivative of  $\phi_{g_{(a,b,c)}}$ 
	at
	$p$ detects all three parameters $a,b,c$, whereas the first-order
	derivative at $p$ detects only $b$. This is precisely why the
	second-order method succeeds in this example while the first-order
	method fails.

	With respect to the local coordinates $(y,z)$ determined by
	$\kappa$, let
	$$
	\mathcal E
	:=
	\bigl(
	e_y,e_z,e_{yy},e_{yz},e_{zz}
	\bigr)
	$$
	be the ordered basis of
	$$
	T_p^{(2)}S^2
	=
	\bigl(J_p^2(S^2,\mathbb R)_0\bigr)^*,
	$$
	where the linear functionals in $\mathcal E$ are defined by
	$$
	\begin{aligned}
		e_y(j_p^2f)
		&=
		\frac{\partial(f\circ\kappa)}{\partial y}(0,0),\\
		e_z(j_p^2f)
		&=
		\frac{\partial(f\circ\kappa)}{\partial z}(0,0),\\
		e_{yy}(j_p^2f)
		&=
		\frac{\partial^2(f\circ\kappa)}{\partial y^2}(0,0),\\
		e_{yz}(j_p^2f)
		&=
		\frac{\partial^2(f\circ\kappa)}{\partial y\,\partial z}(0,0),\\
		e_{zz}(j_p^2f)
		&=
		\frac{\partial^2(f\circ\kappa)}{\partial z^2}(0,0),
	\end{aligned}
	$$
	for every smooth function $f$ defined near $p$ with $f(p)=0$.
	We denote by
	$$
	[L]_{\mathcal E}
	$$
	the $5\times5$  matrix of 
	$L\in\operatorname{GL}\bigl(T_p^{(2)}S^2\bigr)$ with respect to the ordered basis  $\mathcal E$.
	By the definition of the second-order tangent map, a routine
	application of the chain rule yields
	$$
	\begin{aligned}
		D_p^{(2)}\phi_{g_{(a,b,c)}}(e_y)
		&=e_y,\\
		D_p^{(2)}\phi_{g_{(a,b,c)}}(e_z)
		&=be_y+e_z,\\
		D_p^{(2)}\phi_{g_{(a,b,c)}}(e_{yy})
		&=-2ae_y+e_{yy},\\
		D_p^{(2)}\phi_{g_{(a,b,c)}}(e_{yz})
		&=-(c+ab)e_y-ae_z+be_{yy}+e_{yz},\\
		D_p^{(2)}\phi_{g_{(a,b,c)}}(e_{zz})
		&=-2bce_y-2ce_z+b^2e_{yy}+2be_{yz}+e_{zz}.
	\end{aligned}
	$$
	Therefore
	\begin{equation}
	\bigl[\rho_{2,p}(g_{(a,b,c)})\bigr]_{\mathcal E}
	=
	\begin{pmatrix}
		1&b&-2a&-(c+ab)&-2bc\\
		0&1&0&-a&-2c\\
		0&0&1&b&b^2\\
		0&0&0&1&2b\\
		0&0&0&0&1
	\end{pmatrix}.
	\tag{4}
	\end{equation}
	It follows immediately from $(4)$ that
	$$
	\begin{aligned}
		&\rho_{2,p}(g_{(a,b,c)})
		=
		\rho_{2,p}(g_{(a',b',c')})
		\\
		&\quad\Longrightarrow\quad
		\bigl[\rho_{2,p}(g_{(a,b,c)})\bigr]_{\mathcal E}
		=
		\bigl[\rho_{2,p}(g_{(a',b',c')})\bigr]_{\mathcal E}
		\\
		&\quad\Longrightarrow\quad
		a=a',\quad b=b',\quad c=c'
		\\
		&\quad\Longrightarrow\quad
		g_{(a,b,c)}=g_{(a',b',c')}.
	\end{aligned}
	$$
	Thus, $\rho_{2,p}$ is injective.

	Moreover, $\rho_{2,p}(\Gamma)$ is a closed discrete subgroup of
	$\operatorname{GL}\bigl(T_p^{(2)}S^2\bigr)$. Here, discreteness is
	understood with respect to the subspace topology inherited from
	$\operatorname{GL}\bigl(T_p^{(2)}S^2\bigr)$.
	For any $g_{(a,b,c)}\in\Gamma$, it follows from $(4)$ that
	$$
	\begin{aligned}
		a&=-
		\bigl(\bigl[\rho_{2,p}(g_{(a,b,c)})\bigr]_{\mathcal E}\bigr)_{2,4},
		\\
		b&=\bigl(\bigl[\rho_{2,p}(g_{(a,b,c)})\bigr]_{\mathcal E}\bigr)_{1,2},
		\\
		c&=-\frac12
		\bigl(\bigl[\rho_{2,p}(g_{(a,b,c)})\bigr]_{\mathcal E}\bigr)_{2,5}.
	\end{aligned}
	$$
	Thus, for any convergent sequence 
	$\rho_{2,p}(g_{(a_n,b_n,c_n)})$
	in 
	$\operatorname{GL}\bigl(T_p^{(2)}S^2\bigr)$,
	the three integer sequences
	$$
	a_n,
	\quad
	b_n,
	\quad
	c_n
	$$
	converge and are therefore eventually constant.
	Hence every convergent sequence in $\rho_{2,p}(\Gamma)$ is eventually
	constant. Since $\operatorname{GL}\bigl(T_p^{(2)}S^2\bigr)$ is metrizable, it
	follows that
	$\rho_{2,p}(\Gamma)$
	is a closed discrete subgroup of $\operatorname{GL}\bigl(T_p^{(2)}S^2\bigr)$.

	Now, let  $\mu_0$  be the Dirac probability measure concentrated
	at $p$, i.e.,  for every Borel set $E\subseteq S^2$,
	$$
	\mu_0(E)=
	\begin{cases}
		1,&p\in E,\\
		0,&p\notin E.
	\end{cases}
	$$
	Fix any continuous section
	$$	\textsl{g}_2:S^2\longrightarrow\operatorname{Riem}^{(2)}(S^2).
	$$

	Evaluation at $p$ induces a canonical isometric bijection
	$$
	\operatorname{ev}_p:
	L^2_{\textsl{g}_2}
	\bigl(
	S^2,\mu_0,\operatorname{Riem}^{(2)}(S^2)
	\bigr)
	\longrightarrow
	\mathcal P\bigl(T_p^{(2)}S^2\bigr),
	\qquad
	\operatorname{ev}_p([s])=s(p).
	$$
	This map is well-defined since $\mu_0(\{p\})=1$.
	By the definition of the $L^2$-metric,
	we have
	$$
	\begin{aligned}
		d_{L^2}([s],[t])^2
		&=
		\int_{S^2}
		d_{\mathcal P(T_x^{(2)}S^2)}
		\bigl(s(x),t(x)\bigr)^2
		\,d\mu_0(x)
		\\
		&=
		d_{\mathcal P(T_p^{(2)}S^2)}
		\bigl(s(p),t(p)\bigr)^2.
	\end{aligned}
	$$
	Taking square roots gives
	$$
	d_{L^2}([s],[t])
	=
	d_{\mathcal P(T_p^{(2)}S^2)}
	\bigl(s(p),t(p)\bigr).
	$$
	Therefore $\operatorname{ev}_p$ is an isometric embedding and, in
	particular, is injective.
	It remains to prove surjectivity. Let
	$$
	Q\in\mathcal P\bigl(T_p^{(2)}S^2\bigr).
	$$
	Define a measurable section $s_Q$ by
	$$
	s_Q(x)
	=
	\begin{cases}
		Q,&x=p,\\
		\textsl{g}_{2,x},&x\neq p.
	\end{cases}
	$$
	Since
	$$
	\begin{aligned}
		\int_{S^2}
		d_{\mathcal P(T_x^{(2)}S^2)}
		\bigl(s_Q(x),\textsl{g}_{2,x}\bigr)^2
		\,d\mu_0(x)
		&=
		d_{\mathcal P(T_p^{(2)}S^2)}
		\bigl(Q,\textsl{g}_{2,p}\bigr)^2
		\\
		&<\infty,
	\end{aligned}
	$$
	we have
	$$
	[s_Q]\in  L^2_{\textsl{g}_2}
	\bigl(
	S^2,\mu_0,\operatorname{Riem}^{(2)}(S^2)
	\bigr),
	$$
	and
	$$
	\operatorname{ev}_p([s_Q])=Q.
	$$
	Thus $\operatorname{ev}_p$ is surjective and hence an isometric
	bijection.

	Consider the composite map
	$$
	\mathcal{F}
	:=
	\operatorname{ev}_p\circ F_{\mu_0,\textsl{g}_2}^{(2)}\big|_\Gamma
	:
	\Gamma
	\longrightarrow
	\mathcal P\bigl(T_p^{(2)}S^2\bigr).
	$$
	Since every element of $\Gamma$ fixes $p$,  $\mathcal{F}$ is given explicitly by
	\begin{equation}
	\begin{aligned}
		\mathcal{F}(\gamma)
		&=
		\bigl((\phi_\gamma)_*^{(2)}\textsl{g}_2\bigr)_p\\
		&=
		\bigl(D_p^{(2)}\phi_{\gamma^{-1}}\bigr)^*
		\textsl{g}_{2,\phi_{\gamma^{-1}}(p)}\\
		&=
		\bigl(\rho_{2,p}(\gamma)^{-1}\bigr)^*\textsl{g}_{2,p}
	\end{aligned}\tag{5} 
	\end{equation}
	for every $\gamma\in\Gamma$. 
	Since $\operatorname{ev}_p$ is an isometric
	bijection, it suffices to show that
	$\mathcal{F}$ is a coarse embedding.

	For this purpose, we need to show that for any compact subset 
	$\mathcal C\subseteq\mathcal P(T_p^{(2)}S^2)$, 
	$\mathcal{F}^{-1}(\mathcal C)$ is finite. 
	Consider the composition
	$$
	\pi:
	\operatorname{GL}(T_p^{(2)}S^2)
	\xlongrightarrow{\pi_0}
	{\operatorname{GL}(T_p^{(2)}S^2)}\big/
	{\operatorname{O}(T_p^{(2)}S^2,\textsl{g}_{2,p})}
	\xlongrightarrow{\sigma}
	\mathcal P(T_p^{(2)}S^2),
	$$
	where $\pi_0$ is the canonical quotient projection and $\sigma$ is
	the natural homeomorphism given by
	$$
	\sigma\bigl(L\cdot\operatorname{O}(T_p^{(2)}S^2,\textsl{g}_{2,p})\bigr)
	:=(L^{-1})^*\textsl{g}_{2,p}.
	$$
	Since $\operatorname{O}(T_p^{(2)}S^2,\textsl{g}_{2,p})$ is compact,
	$\pi_0$ is proper. 
	Hence $\pi=\sigma\circ\pi_0$ is proper.
	Note that 
	$\mathcal{F}=\pi\circ\rho_{2,p}$.
	Let $\mathcal C\subseteq\mathcal P(T_p^{(2)}S^2)$ be compact.
	Since $\pi$ is proper and $\rho_{2,p}(\Gamma)$ is closed and discrete in $\operatorname{GL}(T_p^{(2)}S^2)$,
	the set
	$$
	\rho_{2,p}(\Gamma)\cap\pi^{-1}(\mathcal C)
	$$
	is finite. Since $\rho_{2,p}$ is injective, it follows that
	$\mathcal{F}^{-1}(\mathcal C)
	=
	\rho_{2,p}^{-1}
	\left(
	\rho_{2,p}(\Gamma)\cap\pi^{-1}(\mathcal C)
	\right)$ 
	is finite.

	We now show that $\mathcal{F}$ is a coarse embedding. 
	
	For the upper control, define
	$$
	C:=
	\max_{\gamma\in\Delta}
	d_{\mathcal P(T_p^{(2)}S^2)}
	\bigl(\textsl{g}_{2,p},\mathcal{F}(\gamma)\bigr)
	<\infty,
	$$
	where $\Delta$ is the finite symmetric generating set.
	Given $\gamma\in\Gamma$, choose a shortest word
	$$
	\gamma=\gamma_1\cdots\gamma_n,
	\qquad
	\gamma_j\in\Delta,
	\qquad
	n=d_\Gamma(1_\Gamma,\gamma),
	$$
	and set $\eta_0:=1_\Gamma$ and
	$\eta_j:=\gamma_1\cdots\gamma_j$ for $1\leq j\leq n$.
	Then
	$$
	\eta_n=\gamma,
	\qquad
	\eta_{j-1}^{-1}\eta_j=\gamma_j
	\quad (1\leq j\leq n).
	$$
	Since $\mathcal{F}(\eta_0)=\mathcal{F}(1_\Gamma)=\textsl{g}_{2,p}$
	and $\eta_n=\gamma$, the triangle inequality gives
	\begin{equation}
	d_{\mathcal P(T_p^{(2)}S^2)}
	\bigl(\textsl{g}_{2,p},\mathcal{F}(\gamma)\bigr)
	\leq
	\sum_{j=1}^{n}
	d_{\mathcal P(T_p^{(2)}S^2)}
	\bigl(\mathcal{F}(\eta_{j-1}),\mathcal{F}(\eta_j)\bigr).
	\tag{6}
	\end{equation}
	Moreover, for each $j=1,\ldots,n$, equation $(5)$ and the fact that
	$\rho_{2,p}$ is a group homomorphism yield 
	\begin{equation}
	\begin{aligned}
		&
		d_{\mathcal P(T_p^{(2)}S^2)}
		\bigl(\mathcal{F}(\eta_{j-1}), \enspace \mathcal{F}(\eta_j)\bigr)
		\\
		&=
		d_{\mathcal P(T_p^{(2)}S^2)}
		\left(
		\bigl(\rho_{2,p}(\eta_{j-1})^{-1}\bigr)^*\textsl{g}_{2,p},\enspace
		\bigl(\rho_{2,p}(\eta_j)^{-1}\bigr)^*\textsl{g}_{2,p}
		\right)
		\\
		&=
		d_{\mathcal P(T_p^{(2)}S^2)}
		\left(
		\textsl{g}_{2,p}, \enspace
		\rho_{2,p}(\eta_{j-1})^*
		\left(
		\bigl(\rho_{2,p}(\eta_j)^{-1}\bigr)^*\textsl{g}_{2,p}
		\right)
		\right)
		\\
		&=
		d_{\mathcal P(T_p^{(2)}S^2)}
		\left(
		\textsl{g}_{2,p}, \enspace
		\bigl(
		\rho_{2,p}(\eta_j)^{-1}\rho_{2,p}(\eta_{j-1})
		\bigr)^*\textsl{g}_{2,p}
		\right)
		\\
		&=
		d_{\mathcal P(T_p^{(2)}S^2)}
		\left(
		\textsl{g}_{2,p}, \enspace
		\bigl(
		\rho_{2,p}(\eta_{j-1}^{-1}\eta_j)^{-1}
		\bigr)^*\textsl{g}_{2,p}
		\right)
		\\
		&=
		d_{\mathcal P(T_p^{(2)}S^2)}
		\bigl(
		\textsl{g}_{2,p}, \enspace
		\mathcal{F}(\eta_{j-1}^{-1}\eta_j)
		\bigr)
		\\
		&=
		d_{\mathcal P(T_p^{(2)}S^2)}
		\bigl(\textsl{g}_{2,p}, \enspace \mathcal{F}(\gamma_j)\bigr)
		\\
		&\leq C.
	\end{aligned}
	\tag{7}
	\end{equation}
	Thus, it follows from  $(6)$ and $(7)$ that 
	$$
	\begin{aligned}
		d_{\mathcal P(T_p^{(2)}S^2)}
		\bigl(\textsl{g}_{2,p},\mathcal{F}(\gamma)\bigr)
		\leq
		C\,d_\Gamma(1_{\Gamma}, \gamma),
	\end{aligned}
	$$
	and hence for arbitrary $\gamma,\eta\in\Gamma$, 
	\begin{equation}
	\begin{aligned}
		&
		d_{\mathcal P(T_p^{(2)}S^2)}
		\bigl(\mathcal{F}(\gamma),\mathcal{F}(\eta)\bigr)
		\\
		&=
		d_{\mathcal P(T_p^{(2)}S^2)}
		\bigl(\textsl{g}_{2,p},\mathcal{F}(\gamma^{-1}\eta)\bigr)
		\\
		&\leq
		C\,d_\Gamma(1_{\Gamma}, \gamma^{-1}\eta)
		\\
		&=
		C\,d_\Gamma(\gamma,\eta).
	\end{aligned}
	\tag{8}
	\end{equation}

	For the lower control, define
	$$
	\rho_-(R)
	:=
	\inf_{\substack{
			\gamma\in\Gamma\\
			d_\Gamma(1_{\Gamma}, \gamma) \geq R
	}}
	d_{\mathcal P(T_p^{(2)}S^2)}
	\bigl(\textsl{g}_{2,p},\mathcal{F}(\gamma)\bigr).
	$$
	We have shown that the preimage  of every compact subset
	of
	$\mathcal P\bigl(T_p^{(2)}S^2\bigr)$
	under $\mathcal{F}$
	is finite. 
	Since
	$\mathcal P(T_p^{(2)}S^2)$ is a finite-dimensional complete
	Riemannian manifold, every closed metric ball in this space is compact.
	It follows that, for every $T\geq0$,
	$$
	K_T
	:=
	\left\{
	\gamma\in\Gamma:
	d_{\mathcal P(T_p^{(2)}S^2)}
	\bigl(\textsl{g}_{2,p},\mathcal{F}(\gamma)\bigr)\leq T
	\right\}
	=
	\mathcal{F}^{-1}\bigl(\overline{B}(\textsl{g}_{2,p},T)\bigr)
	$$
	is finite. Hence
	$$
	R_T
	:=
	1+\max_{\gamma\in K_T}d_\Gamma(1_{\Gamma},\gamma)<\infty.
	$$
	By the definition of $R_T$,
	$$
	d_\Gamma(1_{\Gamma}, \gamma)\geq R_T
	\quad\Longrightarrow\quad
	\gamma\notin K_T
	\quad\Longrightarrow\quad
	d_{\mathcal P(T_p^{(2)}S^2)}
	\bigl(\textsl{g}_{2,p},\mathcal{F}(\gamma)\bigr)>T.
	$$
	Taking the infimum yields
	$$
	\rho_-(R_T)\geq T.
	$$
	Since $\rho_-$ is nondecreasing and $T\geq0$ is arbitrary, it follows
	that
	$$
	\rho_-(R)\longrightarrow\infty
	\quad
	\text{as }R\longrightarrow\infty.
	$$
	For arbitrary $\gamma,\eta\in\Gamma$, we have
	\begin{equation}
	\begin{aligned}
		d_{\mathcal P(T_p^{(2)}S^2)}
		\bigl(\mathcal{F}(\gamma),\mathcal{F}(\eta)\bigr)
		&=
		d_{\mathcal P(T_p^{(2)}S^2)}
		\bigl(\textsl{g}_{2,p},\mathcal{F}(\gamma^{-1}\eta)\bigr)
		\\
		&\geq
		\rho_-\bigl(d_\Gamma(1_{\Gamma}, \gamma^{-1}\eta)\bigr)
		\\
		&=
		\rho_-\bigl(d_\Gamma(\gamma,\eta)\bigr).
	\end{aligned}
	\tag{9}
	\end{equation}

	Combining $(8)$ and $(9)$,  we obtain
	$$
	\rho_-\bigl(d_\Gamma(\gamma,\eta)\bigr)
	\leq
	d_{\mathcal P(T_p^{(2)}S^2)}
	\bigl(\mathcal{F}(\gamma),\mathcal{F}(\eta)\bigr)
	\leq
	C\,d_\Gamma(\gamma,\eta),
	$$
	where $\rho_-(R)\to\infty$. Therefore $\mathcal{F}$ is a coarse embedding.
\end{proof}

\begin{rmk}\label{rem:property-T-second-order}
The same method also yields, for each $n\geq3$, a Property~$(T)$ subgroup 
$\Gamma\cong\mathbb Z^n\rtimes\operatorname{SL}(n,\mathbb Z)$
of $\operatorname{Diff}(S^n)$  
that is second-order $\mu_0$-discrete for some Borel
probability measure $\mu_0$, but not first-order
$\mu$-discrete for any finite Borel measure $\mu$ on $S^n$.
These examples show that higher-order orbit maps detect
a  larger class of subgroups of $\operatorname{Diff}(M)$
than their lower-order counterparts.
In particular, our theorem generalizes the first-order
$\mu$-discreteness criterion of Gong--Wu--Xie--Yu
\cite{GongWuXieYu} for the Novikov conjecture.
\end{rmk}

\end{document}